\documentclass[11pt,a4paper]{article}
\usepackage{amsmath, amssymb}
\usepackage{latexsym, color}
\usepackage{epsfig}

\numberwithin{equation}{section}

\begin{document}

\newtheorem{thm}{Theorem}[section]
\newtheorem{cor}[thm]{Corollary}
\newtheorem{lem}[thm]{Lemma}
\newtheorem{prop}[thm]{Proposition}
\newtheorem{definition}[thm]{Definition}
\newtheorem{rem}[thm]{Remark}
\newtheorem{Ex}[thm]{EXAMPLE}
\def\nm{\noalign{\medskip}}

\bibliographystyle{plain}

%\numberwithin{equation}{section}

%%%%%%%%%%%%%%%%%%%%%% Greek Letter
\newcommand{\Ga}{\alpha}
\newcommand{\Gl}{\lambda}
\newcommand{\Gs}{\sigma}
\newcommand{\GL}{\Lambda}

%%%%%%%%%%%%%%%%%%%%% Bold Letter
\newcommand{\Bc}{{\bf c}}
\newcommand{\Bn}{{\bf n}}
\newcommand{\Bp}{{\bf p}}
\newcommand{\Bt}{{\bf t}}
\newcommand{\Bu}{{\bf u}}
\newcommand{\Bx}{{\bf x}}
\newcommand{\By}{{\bf y}}

\newcommand{\qed}{\hfill \ensuremath{\square}}
\newcommand{\ds}{\displaystyle}
\newcommand{\pf}{\medskip \noindent {\sl Proof}. ~ }
\newcommand{\p}{\partial}
\renewcommand{\a}{\alpha}
\newcommand{\z}{\zeta}
\newcommand{\pd}[2]{\frac {\p #1}{\p #2}}
\newcommand{\norm}[1]{\| #1 \|}
\newcommand{\dbar}{\overline \p}
\newcommand{\eqnref}[1]{(\ref {#1})}
\newcommand{\na}{\nabla}
\newcommand{\Om}{\Omega}
\newcommand{\ep}{\epsilon}
\newcommand{\tmu}{\widetilde \mu}
\newcommand{\vep}{\varepsilon}
\newcommand{\tlambda}{\widetilde \lambda}
\newcommand{\tnu}{\widetilde \nu}
\newcommand{\vp}{\varphi}
\newcommand{\RR}{\mathbb{R}}
\newcommand{\CC}{\mathbb{C}}
\newcommand{\NN}{\mathbb{N}}
\renewcommand{\div}{\mbox{div}~}

\newcommand{\la}{\langle}
\newcommand{\ra}{\rangle}
\newcommand{\Scal}{\mathcal{S}}
\newcommand{\Lcal}{\mathcal{L}}
\newcommand{\Kcal}{\mathcal{K}}
\newcommand{\Dcal}{\mathcal{D}}
\newcommand{\tScal}{\widetilde{\mathcal{S}}}
\newcommand{\tKcal}{\widetilde{\mathcal{K}}}
\newcommand{\Pcal}{\mathcal{P}}
\newcommand{\Qcal}{\mathcal{Q}}
\newcommand{\id}{\mbox{Id}}

%%%%%%%%%%
\newcommand{\be}{\begin{equation}}
\newcommand{\ee}{\end{equation}}

\title{Asymptotics and Computation of the Solution to the Conductivity Equation in the Presence of Adjacent Inclusions with Extreme Conductivities\thanks{\footnotesize This work was supported by Ministry of Education, Sciences and Technology of Korea through NRF grants No. 2009-0090250 (HK and ML) and 2009-0070442 (ML), and by Hankuk University of Foreign Studies Research Fund of 2011 (KY).}}

\author{Hyeonbae Kang\thanks{\footnotesize Department of Mathematics, Inha University, Incheon 402-751, Republic of Korea (hbkang@inha.ac.kr)} \and Mikyoung Lim\thanks{Department of Mathematics, Korea Advanced Institute of Science and Technology, Yuseong-gu, Daejeon 305-701, Republic of Korea (mklim@kaist.ac.kr)} \and KiHyun Yun\thanks{\footnotesize Department of Mathematics,
Hankuk University of Foreign Studies,  Youngin-si, Gyeonggi-do 449-791, Republic of Korea (gundam@hufs.ac.kr)}}

\date{}

\maketitle

\begin{abstract}
When inclusions with extreme conductivity (insulator or perfect conductor) are closely located, the gradient of the solution to the conductivity equation can be arbitrarily large. And computation of the gradient is extremely challenging due to its nature of blow-up in a narrow region in between inclusions. In this paper we characterize explicitly the singular term of the solution when two circular inclusions with extreme conductivities are adjacent. Moreover, we show through numerical computations that the characterization of the singular term can be used efficiently for computation of the gradient in the presence adjacent inclusions.
\end{abstract}

\noindent {\footnotesize {\bf Mathematics subject classification
(MSC2000): 35J25, 73C40} }


\section{Introduction}

Frequently in composites which consist of inclusions and background (the matrix), the inclusions are
closely spaced, and it is quite important from a practical point of view to know whether the gradient of the potential can be arbitrarily large as the inclusions get closer to each other.  The gradient of the potential represents the stress in anti-plane elasticity and the electric field in the conductivity problem; see \cite{bab}. It is known that the gradient of the potential may blow up
as the distance between the inclusions goes to zero and their material parameters (conductivities or stiffness)
degenerate.

Suppose that $B_1$ and $B_2$ are inclusions whose conductivity is $k$. We suppose that the conductivity of the background is $1$ ($k \neq 1$). Let $\ep$ be the distance between $B_1$ and $B_2$ and assume that $\ep$ is small.
The problem is to estimate $|\nabla u|$, where $u$ is the electrical potential, in terms of $\ep$ when $\ep$ tends to $0$.

There have been important works on this problem. If $k$ stays away from $0$ and $\infty$, {\it i.e.}, $c_1 < k < c_2$ for some positive constants $c_1$ and $c_2$, then it was proved by Bonnetier-Vogelius \cite{BV} and Li-Vogelius \cite{LV} that $|\nabla u|$ remains bounded regardless of $\ep$. This result was extended to elliptic system by Li-Nirenberg \cite{LN}. It is worth emphasizing that the results in \cite{LV, LN} are not only for two inclusions case but also for the case of arbitrary number of inclusions.

On the other hand, if $k$ is either $0$ (insulating) or $\infty$ (perfectly conducting), then $\nabla u$ may blow up as $\ep$ tends to $0$. For two identical perfectly conducting circular inclusions it was shown in \cite{BC} (see also \cite{mar} and \cite{K}) that the
gradient in general becomes unbounded as $\epsilon$ approaches zero and the blow-up rate is
$\epsilon^{-1/2}$. In \cite{AKL, AKLLL}, a lower bound and an upper bound for the gradient has been obtained. These bounds are valid for all $k$ including extreme values ($k=0$ and $k=\infty$) and provide the precise dependence of $\nabla u$ on $\ep$, $k$ and radii of disks. The blow-up of the gradient may or may not occur depending on the background potential. In \cite{AKLLZ}, Ammari {\it et al} characterize those background potential which actually make the gradient blow up. In \cite{Y, Y2}, Yun showed that the blow-up rate is $\ep^{-1/2}$ for perfectly conducting and insulated inclusions of arbitrary shape in two dimensions. In three dimensions, Bao {\it et al} \cite{BLY} proved that the blow-up rate for the perfectly conducting inclusions is $|\ep \log \ep|^{-1}$ and extended the result to the case of multiple inclusions \cite{BLY2}. Lim-Yun \cite{LY} also found the same blow-up rate when inclusions are spheres. Their estimates explicitly reveal the dependence on the radii of the sphere. They also showed in \cite{LY2} that if there is a small bump in between two inclusions in two dimensions, then the magnitude of the blow-up gets larger.

The purpose of this paper is to characterize the singular term of the solution, {\it i.e.}, to establish an asymptotic formula for the blow-up of the gradient when two circular inclusions get closer. We find the decomposition of the solution $u$ to the conductivity equation as
 \be\label{decompo}
 u=g+b
 \ee
where $\nabla g$ may blow up at the rate of $\ep^{-1/2}$ while $\nabla b$ stays bounded regardless of $\ep$, when $B_1$ and $B_2$ are disks and $k$ is either $\infty$ or $0$. We actually obtain an explicit formula for the term $g$ which gives a precise description of singular behavior of $\nabla u$.

The characterization of the singular term of the solution finds a very good application in the computation of electrical fields. Computation of the electrical field in the presence of closely located inclusions with extreme ($0$ or $\infty$) conductivities is known to be a extremely difficult problem because of the the blow-up phenomenon in a very narrow region between inclusions. Since the the gradient of the solution is arbitrarily large, we need very fine mesh to catch the large gradient in a narrow region. The results of this paper constitute a significant step toward overcoming this difficulty since the singular term $g$ is explicit and computation of $b$ requires only regular meshes. We present efficient methods to use the decomposition for the computation of the solution and some results of numerical computation using them. Numerical examples of this paper show that these methods work pretty well.

This paper is organized as follows. In the next section we derive the decomposition \eqnref{decompo} for the perfect conductors in the free space. In section 3, we deal with the same problem in bounded domains. Section 4 is for the insulators. New numerical methods and results of computation are presented in the last section.

The result of this paper can be extended to perfect conductors of spherical shape in three dimensions. This result will be presented in a forthcoming paper.

%%%%%%%%%%%%%%%%%%%%%%%%%%%%%%%%%%%%%%%%%%%%%%%%%%%%%%
\section{Free space problem-perfectly conducting case}
%%%%%%%%%%%%%%%%%%%%%%%%%%%%%%%%%%%%%%%%%%%%%%%%%%%%%%

Let $B_j=B(\Bc_j, r_j)$, $j=1,2$, be the disk centered at $\Bc_j$ and of radius $r_j$, and
 \be
 \sigma =
 \left\{
 \begin{aligned}
 & 1 \quad \mbox{on } \RR^2 \setminus (B_1 \cup B_2), \\
 & k \quad \mbox{on } B_1 \cup B_2,
 \end{aligned}
 \right .
 \ee
which represents the conductivity distribution: the conductivity of the inclusions is $k$ ($k \neq 1$) and that of the background is 1. The equation we consider is
 \be\label{condeqn}
 \nabla \cdot \Gs \nabla u =0 \quad\mbox{in } \RR^2,
 \ee
which may be viewed as the conductivity equation or anti-plane elasticity equation. A condition at the infinity is prescribed by
 \be\label{infcond}
 u(\Bx)-H(\Bx) = O(|\Bx|^{-1}) \quad \mbox{as } |\Bx| \to \infty,
 \ee
where $H$ is an entire harmonic function and represents the background potential.

If $k=\infty$, the equation \eqnref{condeqn} with the condition \eqnref{infcond} is understood as the following problem:
\begin{equation} \label{eq:main}
\left\{
\begin{array}{ll}
\ds\Delta u  = 0\quad&\mbox{in }{\mathbb{R}^2 \setminus \overline{B_1 \cup B_2}},\\
\ds u |_{\partial B_j} = \lambda_j \ (\mbox{constant}) \quad& j=1,2,\\
\ds u(\Bx)- H(\Bx) = O(|\Bx|^{-1})\quad&\mbox{as } |\Bx| \rightarrow \infty.
\end{array}
\right.
\end{equation}
The constants $\lambda_j$ can be determined by the additional requirements
 \be
 \int_{\p B_j} \frac{\p u}{\p \nu} \Big|_{+} ~ds  = 0, \quad j=1,2,
 \ee
where $\nu$ is the outward unit normal vector of ${\mathbb{R}^2 \setminus \overline{B_1 \cup B_2}}$,
{\it i.e.}, directed inward of $B_i$. Here and throughout this paper, the notations $|_{+}$ and $|_{-}$ are for limits from outside and inside inclusions, respectively.

Let $R_j$, $j=1,2$, be the reflection
with respect to $\p B_j$, {\em i.e.},
\be\label{ref_B}
R_j (\Bx): = \frac{r_j^2 (\Bx-\Bc_j)}{|\Bx-\Bc_j|^2}+\Bc_j, \quad j=1,2.
\ee
It is easy to see that the combined reflections $R_1R_2$ and
$R_2R_1$ have unique fixed points, say $\Bp_1$ and $\Bp_2$, respectively. Let
 \be\label{hBx}
 h(\Bx):=\frac 1 {2\pi} \left( \log |\Bx-\Bp_1|- \log |\Bx-\Bp_2|\right).
 \ee
The function $h$, which was first found in \cite{LY}, has a special property: it is the solution to
\be\label{eqn:h}
\begin{cases}
\ds\Delta h=0 \quad& \mbox{in } \RR^2 \setminus \overline{B_1 \cup B_2}, \\
\ds h|_{\p B_j}=  C_j \ (\mbox{constant}) \quad& j=1,2,\\
\nm
\ds \int_{\p B_j} \pd{h}{\nu} ~ds=(-1)^j, \quad& j=1,2,\\
\ds h(\Bx)=O(|\Bx|^{-1}) \quad&\mbox{as } |\Bx| \to \infty.
\end{cases}
\ee
The following formula was proved in \cite{Y,Y2}: let $\Gl_1$ and $\Gl_2$ be constants appearing in \eqnref{eq:main}, then
 \be\label{potdiff}
 \Gl_2-\Gl_1 = \int_{\partial B_1 \cup \partial B_2 } H \partial_{\nu} h ~ds  = H(\textbf{p}_2) - H(\textbf{p}_1).
 \ee

The following is the first main theorem of this paper.

\begin{thm} \label{Stiff}
Let $\Bn$ be the unit vector in the direction of $\Bp_2-\Bp_1$ and let $\Bp$ be the middle point of the shortest line segment connecting $\p B_1$ and $\p B_2$. For a harmonic function $H$ in $\RR^2$, let $u$ be the solution to \eqref{eq:main}. Then, the solution $u$ can be expressed as follows:
 \be\label{uahb}
 u(\Bx)= a h(\Bx) + b(\Bx)
 \ee
where
 \be \label{intensity}
 a= \frac {4\pi r_1 r_2}{r_1 + r_2} (\Bn \cdot \nabla H) (\Bp)
 \ee
and for any bounded set $\Om$ containing $\overline{B}_1$ and $\overline{B}_2$ there is a constant $C$ independent of $\ep$ such that
 \be
 \norm{\nabla b}_{L^{\infty} (\Om \setminus ( B_1 \cup B_2))} \leq C .
 \ee
\end{thm}

The asymptotic formula $\nabla u$ as $\ep \to 0$ is then given by
 \be
 \nabla u(\Bx)= \frac {2 r_1 r_2}{r_1 + r_2} (\Bn \cdot \nabla H) (\Bp)  \left( \frac {\mathbf{x}-\mathbf{p}_1} { \bigr|\mathbf{x}-\mathbf{p}_1\bigr|^2}- \frac {\mathbf{x}-\mathbf{p}_2} {\bigr|\mathbf{x}-\mathbf{p}_2\bigr|^2}\right) + O(1). \label{asymp}
 \ee

Let us make a few remarks on Theorem \ref{Stiff} before proving it. It is shown in \cite{Y,Y2} that the fixed points $\Bp_1$ and $ \Bp_2$ are given by
 $$
 \Bp_1 = \Bigr(-\sqrt{2}\sqrt{\frac{r_1r_2}{r_1+r_2}}\sqrt\epsilon +O(\epsilon)
 ,0\Bigr) \quad \mbox{and} \quad \Bp_2 =
 \Bigr(\sqrt{2}\sqrt{\frac{r_1r_2}{r_1+r_2}}\sqrt\epsilon +O(\epsilon)
 ,0\Bigr)
 $$
if $\mathbf{c}_1=(-r_1 - \frac{\epsilon}{2}, 0)$ and $\mathbf{c}_2=(r_2+ \frac{\epsilon}{2},0)$.
If $r_1$ and $r_2$ are bounded below by a positive constant $r_0$, then there are positive constants $C_1$ and $C_2$ depending only on $r_0$ such that
 \be\label{nablah-1}
 C_1\sqrt {\frac {r_1 + r_2}{r_1 r_2}} \frac{1}{\sqrt{\ep}}
\le |\nabla h (\Bx)| \le C_2\sqrt {\frac {r_1 + r_2}{r_1 r_2}} \frac{1}{\sqrt{\ep}}
 \ee
for all $\Bx$ on the the shortest line segment connecting $\p B_1$ and $\p B_2$, see \cite{LY2}. Thus the blow-up rate of $|\nabla u|$ is $\ep^{-1/2}$. It is also proved in the same paper that
 \be\label{nablah2}
 |\nabla h (\Bx)| \le C\sqrt {\frac {r_1 + r_2}{r_1 r_2}}  \frac{1}{\sqrt{\ep}}
 \ee
for all $\Bx$. Thus, an optimal bound for $\nabla u$ in a bounded domain can be obtained from \eqref{nablah-1} and  \eqref{nablah2} in terms of $r_1$, $r_2$, $\epsilon$ and  $(\Bn \cdot \nabla H) (\Bp)$. In view of the formula \eqnref{intensity} of $a$ (we call it the stress intensity factor), the blow-up does not occur if $(\Bn \cdot \nabla H) (\Bp)=0$. This fact was already found in \cite{AKLLZ}. One can also show that if $r_1$ and $r_2 $ are $ O(\ep)$, then there is a constant $C$ independent of $\ep$ such that
 \be
 |\nabla h (\Bx)| \le \frac{C}   {\ep}
 \ee
for all $\Bx$. Thus \eqref{asymp} means that in this case, no blow-up occurs: $\nabla u$ stays bounded. This finding is in agreement with that in \cite{AKLLL}.

We first prove the following proposition by modifying an argument of Bao {\it et al} \cite{BLY}.
\begin{prop}\label{prop:boundedness}
Let
 \be\label{decom}
 b(\Bx) = u(\Bx)- \left(\frac { u |_{\partial B_2} -u |_{\partial B_1} }{h |_{\partial B_2} - h |_{\partial B_1} }\right) h(\Bx), \quad \Bx \in \RR^2 \setminus (B_1 \cup B_2).
 \ee
For any bounded set $\Om_1$ containing $\overline{B}_1$ and $\overline{B}_2$ and $\Om_2$ containing $\overline{\Om}_1$, there is a constant $C$ independent of $\ep$ such that
 \be
 \norm{b}_{\mathcal{C}^1(\Om_1 \setminus {(B_1 \cup B_2)})} \leq C \norm{H}_{L^{\infty}(\Om_2)  }.
 \ee
\end{prop}

\pf
It can be easily seen that $b$ is bounded. In fact, since $b$ is harmonic in $\RR^2 \setminus \overline{B_1 \cup B_2}$ and
 $$
 b|_{\p B_2} - b|_{\p B_1} =0,
 $$
we infer from the result in \cite{ADKL} (see also \cite{BLY}) that $b$ is bounded in $\Om_1$.

Since $h(\Bx) \to 0$ as $|\Bx| \to \infty$, by the maximum principle, $h$ attains its maximum and minimum on $\p B_2$ and $\p B_1$, respectively.  Thus, we have
 $$
 \norm {h}_{L^{\infty}(\mathbb{R}^2 \setminus (B_1 \cup B_2))}\leq h |_{\partial B_2} - h |_{\partial B_1},
 $$
and hence
 \be\label{guuH}
 \norm {b}_{L^{\infty}(\Om_2 \setminus  {(B_1 \cup B_2)})  } \leq  \norm {u}_{L^{\infty}(\Om_2 \setminus (B_1 \cup B_2))  } +\left| u |_{\partial B_2} - u |_{\partial B_1} \right|.
 \ee

Let $\Gl_j = u|_{\p B_j}$, $j=1,2$, and assume that $\Gl_2 \ge \Gl_1$ without loss of generality. Since $(u-H)(\Bx) \to 0$ as $|\Bx| \to \infty$, the maximum and minimum of $u-H$ occur on $\p B_1 \cup \p B_2$. So, we have
 $$
 (u-H)(\Bx) \le \max_{\p B_1 \cup \p B_2 }(u-H) \le \Gl_2 +  \norm{H}_{L^{\infty}(\Om_2)},
 $$
and
 $$
 \Gl_1 -  \norm {H}_{L^{\infty}(\Om_2)}  \le \min_{\p B_1 \cup \p B_2 }(u-H) \leq (u -H)(\Bx)
 $$
for all $\Bx \in \RR^2 \setminus ({B_1 \cup B_2})$. Since $\min_{\mathbb{R}^2 \setminus (B_1 \cup B_2)} (u-H) \le 0$,
we have
 $$
 \Gl_1 \le  \norm {H}_{L^{\infty}(\Om_2)}.
 $$
Therefore, we have
 \begin{align*}
 \norm {u-H}_{L^{\infty}(\Om_2 \setminus (B_1 \cup B_2))}  & \le \Gl_2 +  \norm{H}_{L^{\infty}(\Om_2)} \\
 & \le \Gl_2 - \Gl_1 + 2\norm{H}_{L^{\infty}(\Om_2)} \le 4\norm {H}_{L^{\infty}(\Om_2)},
 \end{align*}
where the last inequality comes from \eqnref{potdiff}. Thus
 \begin{align*}
 \norm{u}_{L^{\infty}(\Om_2 \setminus {(B_1 \cup B_2)})} \le \norm{H}_{L^{\infty}(\Om_2)} + \norm{u-H}_{L^{\infty}(\Om_2 \setminus {(B_1 \cup B_2)})} \le 5 \norm{H}_{L^{\infty}(\Om_2)}.
 \end{align*}
It then follows from \eqnref{guuH} and \eqnref{potdiff} that
 \be
 \norm {b}_{L^{\infty}(\Om_2 \setminus (B_1 \cup B_2))  } \leq 7 \norm {H}_{L^{\infty}(\Om_2)}.\label {b}
 \ee

We now show that
 \be\label{gradb}
 \norm {\nabla b}_{L^{\infty}(\Om_1 \setminus (B_1 \cup B_2))  } \leq C \norm {H}_{L^{\infty}(\Om_2)}.
 \ee
For that purpose, we define the harmonic functions $G_+$ and $G_{-}$ as follows:
\begin{equation*}
\quad \left\{
\begin{array}{ll}
\ds\Delta G_{\pm}  = 0,\quad&\mbox{in }{\Om_2 \setminus \overline{B_1 \cup B_2}},\\
\ds G_{\pm}= \pm\norm {b}_{L^{\infty} (\Om_2 \setminus (B_1 \cup B_2))}\quad&\mbox{on } \p \Om_2,\\
\ds G_{\pm} = b,\quad&\mbox{on } \partial B_1 \cup \partial B_2.
\end{array}
\right.
\end{equation*}
Then, $\pm (G_{\pm}-b) \geq 0$ in $\Om_2 \setminus \overline{B_1 \cup B_2}$ and $G_{\pm}-b =0 $ on $\p B_1 \cup \p B_2$. By Hopf's Lemma, we have
\be
\partial_{\nu} G_+ \leq \partial_{\nu} b \leq \partial_{\nu} G_- \quad\mbox{on}~ \partial B_1 \cup \partial B_2. \label{c}
\ee
We introduce more harmonic functions $G_{+1}$, $G_{+2}$, $G_{-1}$ and $G_{-2}$ defined as follows: for $i=1,~2$,
\begin{equation} \notag
\quad \left\{
\begin{array}{ll}
\ds\Delta G_{\pm i }  = 0,\quad&\mbox{in }{\Om_2 \setminus \overline{B_i}},\\
\ds G_{\pm i }= G_{\pm}= \pm\norm {b}_{L^{\infty} (\Om_2 \setminus (B_1 \cup B_2))}\quad&\mbox{on } \p \Om_2,\\
\ds G_{\pm i} = G_{\pm} = b,\quad&\mbox{on } \p B_i.
\end{array}
\right.
\end{equation}
Since $b|_{\p B_1} = b|_{\p B_2}=\mbox{constant}$, we have
 $$
 G_{+i} (\Bx) \ge b|_{\p B_1 \cup \p B_2}.
 $$
In particular,
 $$
 G_{+i} (\Bx) \ge b(\Bx) = G_+(\Bx) \quad \mbox{on } \p B_1 \cup \p B_2.
 $$
Since $G_{+ i}|_{\p \Om_2}= G_+ |_{\p\Om_2}$, we have
 $$
 G_{+i}- G_+ \geq 0 \quad \mbox{in } \Om_2 \setminus {B_1 \cup B_2}.
 $$
Since $G_{+i}- G_+ =0$ on $\p B_i$, it follows from the Hopf's Lemma that
 \be\label{plus}
 \p_{\nu} G_{+i} \leq \p_{\nu} G_{+} \quad\mbox{on } \p B_i, \ \ i=1,2.
 \ee
Similarly, one can show that
 \be\label{minus}
 \p_{\nu} G_{-i} \geq \p_{\nu} G_{-} \quad \mbox{on } \p B_i, \ \ i=1,2.
 \ee

Note that $ {G_{\pm 1}} /{\norm {b}_{L^{\infty} (\Om_2 \setminus (B_1 \cup B_2))}}$ is a harmonic function in $\Om_2 \setminus \overline{B_1}$ which is $\pm 1 $ on $\p \Omega_2$ and also has a constant value between $-1$ and $1$ on $\p B
_1$, and that $\mbox{dist}({B_1, \p \Om_2}) > c_0$. Thus, $ {G_{\pm 1}} /{\norm {b}_{L^{\infty} (\Om_2 \setminus (B_1 \cup B_2))}}$ can be extended as a harmonic function into $\Om_2 \setminus\overline { B ( \Bc_1, \tilde {r})}$ where $\Bc_1$ is the center of $B_1$ and $\tilde {r}$ is strictly less than the radius of $B_1$ independently of $\epsilon$. Then, we have from interior regularity estimates for elliptic equations and \eqnref{b} that
 $$
 \| \p_{\nu} G_{\pm 1} \|_{L^\infty(\p B_1)} \le C \| b \|_{L^{\infty} (\Om_2 \setminus (B_1 \cup B_2))} \le 7C \norm {H}_{L^{\infty}(\Om_2)}.
 $$
It then follows from \eqnref{c}, \eqnref{plus} and \eqnref{minus} that
 $$
 \| \p_{\nu} b \|_{L^\infty(\p B_1)} \le C \norm {H}_{L^{\infty}(\Om_2)}
 $$
for some constant $C$ independent of $\ep$. Similarly one can show that
 $$
 \| \p_{\nu} b \|_{L^\infty(\p B_2)} \le C \norm {H}_{L^{\infty}(\Om_2)}.
 $$
Since $b$ is constant on $\p B_1$ and $\p B_2$, we get
 \be
 \| \nabla b \|_{L^\infty(\p B_1 \cup \p B_2)} \le C \norm {H}_{L^{\infty}(\Om_2)}.
 \ee
The standard interior regularity estimate for harmonic functions shows that
 $$
 \norm{\nabla b}_{L^{\infty}(\p \Om_1)  } \leq  C \norm{H}_{L^{\infty}(\Om_2)  }.
 $$
The maximum principle now yields \eqnref{gradb}, and the proof is complete.
\qed

We are now ready to prove Theorem \ref{Stiff}.

\par
{\medskip \noindent {\sl Proof of Theorem \ref{Stiff}}. }  After translation and rotation if necessary, we may assume that $\mathbf{c}_1=(-r_1 - \frac{\epsilon}{2}, 0)$ and $\mathbf{c}_2=(r_2+ \frac{\epsilon}{2},0)$. Then $\Bp=(0,0)$ and $\Bn= (1,0)$. It is proved in \cite{Y,Y2} that
 $$
 \Bp_1 =
 \Bigr(-\sqrt{2}\sqrt{\frac{r_1r_2}{r_1+r_2}}\sqrt\epsilon +O(\epsilon)
 ,0\Bigr) \quad \mbox{and} \quad \Bp_2 =
 \Bigr(\sqrt{2}\sqrt{\frac{r_1r_2}{r_1+r_2}}\sqrt\epsilon +O(\epsilon)
 ,0\Bigr).
 $$
Therefore, we get from \eqnref{potdiff}
 \be
 u |_{\partial B_2} -u |_{\partial B_1} = 2 \sqrt 2 \partial_{x_1} H (0,0) \sqrt {\frac {r_1 r_2}{r_1 + r_2}} \sqrt {\epsilon} + O(\epsilon)
 \ee
as $\ep \to 0$. On the other hand, one can see that
 $$
 h|_{\p B_2} - h|_{\p B_1} = \frac {1}{\sqrt 2 \pi} \sqrt{\frac {r_1+r_2} {r_1r_2}}\sqrt\epsilon + O(\epsilon).
 $$
Therefore, we get from \eqnref{decom}
\begin{align*}
u(\Bx) & = \frac { 2 \sqrt 2 ~\partial_{x_1} H (\Bp) \sqrt {\frac {r_1 r_2}{r_1 + r_2}} \sqrt {\epsilon} + O(\epsilon)}{\frac {1}{\sqrt 2 \pi} \sqrt{\frac {r_1+r_2} {r_1r_2}}\sqrt\epsilon + O(\epsilon)} h(\Bx) + b(\Bx) \\
&= \frac {4 \pi r_1 r_2}{r_1 + r_2} \partial_{x_1} H (\Bp) h(\Bx) + O(\sqrt {\epsilon} ) h(\Bx)  + b(\Bx).
\end{align*}
Note that the gradient of $O(\sqrt {\epsilon} ) h(\Bx)$ term is bounded because of \eqnref{nablah2} and so is $b(\Bx)$ by Proposition \ref{prop:boundedness}. Thus we obtain
\eqnref{uahb} by setting $O(\sqrt {\epsilon} ) h(\Bx)  + b(\Bx)$ to be the new $b(\Bx)$. This completes the proof.
\qed

%%%%%%%%%%%%%%%%%%%%%%%%%%%%%%%%%%%%%%%%%%%%%%%%%%%%%%%%%%%
\section{Boundary value problem-perfectly conducting case}
%%%%%%%%%%%%%%%%%%%%%%%%%%%%%%%%%%%%%%%%%%%%%%%%%%%%%%%%%%%

Let $\Om$ be a bounded domain with $\mathcal{C}^2$-boundary containing two circular perfectly conducting inclusions $B_j=B(\Bc_j, r_j)$, $j=1, 2$. We assume that the inclusions are away from $\p\Om$, namely, there is a constant $c_0$ such that
 \be\label{away}
 \mbox{dist}(B_j, \p\Om) \ge c_0, \quad j=1, 2.
 \ee
%We also assume that the radii of inclusions are not small, {\it i.e.}, there is $r_0>0$ such that
% \be\label{rzero}
% r_j \ge r_0, \quad j=1, \ldots, N.
% \ee

We consider the following boundary value problem:
 \be\label{bvp}
 \quad \left\{
 \begin{array}{ll}
 \ds  \Delta u=0 \quad \mbox{in } \Om \setminus \overline{B_1 \cup B_2},\\
 \ds \pd{u}{\nu} \Big|_{\p\Om}=g, \\
 \ds u= \mbox{constant on } \p B_j, \ j=1, 2, \\
 \nm
 \ds \int_{\p B_j} \pd{u}{\nu} \,ds =0 , \ j=1,2.
 \end{array}
 \right.
 \ee
Here $g \in L^2_0(\Om)$ ($0$ indicates that $\int_{\p\Om} g=0$) and we impose the condition that $\int_{\p\Om} u =0$ for the uniqueness of the solution.

In this section we derive an asymptotic formula similar to \eqnref{uahb} for the problem \eqnref{bvp}. Here we only consider the Neumann problem. But the same arguments work equally well for the Dirichlet problem.

Let $\Lambda_\Om: L^2_0 (\p\Om) \to H^1(\p\Om)$ be the Neumann to Dirichlet (NtD) map, {\it i.e.},
 \be
 \GL_\Om [g] := u|_{\p\Om},
 \ee
where $u$ is the solution to \eqnref{bvp}. Because of the assumption \eqnref{away}, we have
 \be\label{GLg}
 \| \GL_\Om [g] \|_{H^1(\p\Om)} \le C \| g \|_{L^2(\p\Om)}
 \ee
for all $g \in L^2_0(\p\Om)$ for some constant $C$ independent of $\ep$. See for example \cite[Theorem 2.2]{FKS99}

For a bounded domain $B$ with $\mathcal{C}^2$ boundary let $\Scal_B$ and $\Dcal_B$ denote the single and double layer potentials on $B$:
 \begin{align*}
 \Scal_B [\vp] (\Bx) &= \frac{1}{2\pi} \int_{\p B} \ln |\Bx-\By| \vp(\By) \, ds(\By), \quad \Bx \in \RR^2, \\
 \Dcal_B [\vp] (\Bx) &= -\frac{1}{2\pi} \int_{\p B} \frac{\langle \Bx-\By, \nu(\By) \rangle}{|\Bx-\By|^2} \vp(\By) \, ds(\By), \quad \Bx \in \RR^2 \setminus \p B.
 \end{align*}
We note $\Scal_B$ maps, as an operator defined on $\p B$, $\mathcal{C}^{0,\Ga}(\p B)$ into $\mathcal{C}^{1,\Ga}(\p B)$ if $\alpha >0$. Thus if $\vp \in \mathcal{C}^{0,\Ga}(\p B)$, then $\Scal_B[\vp]$ belongs to $\mathcal{C}^{1,\Ga}(\overline B)$ and $\mathcal{C}^{1,\Ga}(\RR^2 \setminus B)$. The single layer potential enjoys the following jump relation
 \be\label{singlejump}
 \pd{(\Scal_{B} [\vp] )}{\nu} \Big|_+ - \pd{(\Scal_{B} [\vp] )}{\nu} \Big|_- = \vp \quad \mbox{on } \p B.
 \ee

It is known that there are harmonic functions $H$ and a pair of potentials $(\vp_1, \vp_2) \in \mathcal{C}^{0,\alpha}_0(\p B_1) \times \mathcal{C}^{0,\alpha}_0(\p B_2)$ ($0$ indicates that the integral of $\vp_j$ over $\p B_j$ is zero) for some $\alpha>0$ such that the solution $u$ to \eqnref{uahb} is represented by
 \be\label{uHS}
 u(\Bx)= H(\Bx) + \Scal_{B_1}[\vp_1](\Bx) + \Scal_{B_2}[\vp_2](\Bx), \quad \Bx \in \Om \setminus (B_1 \cup B_2).
 \ee
In fact, $H$ is given by
 \be\label{defH}
 H(\Bx) = -\Scal_{\Om}[g](\Bx) + \Dcal_\Om [\GL_\Om [g]](\Bx), \quad \Bx \in \Om,
 \ee
and $(\vp_1, \vp_2)$ is the unique solution to
 \be\label{varphi:eqn}
 \left\{ \begin{array}{l}
 \ds \lambda \varphi_1 +\pd{(\Scal_{B_2} [\varphi_2] )}{\nu^{(1)}}
 = -\pd{H}{\nu^{(1)}} \quad\mbox{on } \p B_1, \\
 \nm
 \ds  \pd{(\Scal_{B_1} [\varphi_1] )}{\nu^{(2)}} + \lambda \varphi_2 = -\pd{H}{\nu^{(2)}}
  \quad\mbox{on } \p B_2,
 \end{array}  \right.
 \ee
where $\lambda = \frac{1}{2}$. Here $\nu^{(j)}$ denotes the normal vector to $\nu^{(j)}$, $j=1,2$. See \cite{KS96, KS2000} (also \cite{AKL, book2}).

Let $\Om_1$ and $\Om_2$ subdomains of $\Om$ such that $\overline{\Om}_1 \subset \Om_2$ and $\overline{\Om}_2 \subset \Om$. We further assume that $B_1$ and $B_2$ are still away from $\p \Om_1$, {\it i.e.},
 \be\label{away2}
 \mbox{dist}(B_1 \cup B_2, \p\Om_1) >c_1,
 \ee
for some $c_1>0$. By the Runge approximation, there is a sequence of harmonic functions $H_n$ in $\RR^2$ such that $H_n \to H$ as $n \to \infty$ in $L^\infty (\Om_2)$. For each $n$, let $u_n$ be the solution to \eqnref{eq:main} with $H$ replaced with $H_n$. Then $u_n$ can be represented as
 \be
 u_n (\Bx)= H_n (\Bx) + \Scal_{B_1}[\vp_1^{(n)}](\Bx) + \Scal_{B_2}[\vp_2^{(n)}](\Bx), \quad \Bx \in \RR^2 \setminus (B_1 \cup B_2),
 \ee
where $(\vp_1^{(n)}, \vp_2^{(n)})$ is the unique solution to \eqnref{varphi:eqn} with $H$ replaced with $H_n$.
Since $\pd{H_n}{\nu^{(j)}} \to \pd{H}{\nu^{(j)}}$ as $n \to \infty$ in $\mathcal{C}^{0,\alpha}(\p B_j)$ for $j=1$ and $2$, we infer from the linearity of the integral equation \eqnref{varphi:eqn} that $(\vp_1^{(n)}, \vp_2^{(n)}) \to (\vp_1, \vp_2)$ as $n \to \infty$ in $\mathcal{C}^{0,\alpha}(\p B_1) \times \mathcal{C}^{0,\alpha}(\p B_2)$. It means that $u_n \to u$ in $\mathcal{C}^{1,\alpha}(\Om_1 \setminus (B_1 \cup B_2))$. We thus get the following theorem from Theorem \ref{Stiff}.

\begin{thm} \label{Stiff-bvp}
Let $u$ be the solution to \eqref{bvp} and $H$ be the function defined by \eqref{defH}. Then, the solution $u$ can be expressed as follows:
 \be\label{uahb2}
 u(\Bx)= a h(\Bx) + b(\Bx), \quad \Bx \in \Om \setminus \overline{B_1 \cup B_2}
 \ee
where
 \be \label{intensity2}
 a= \frac {4\pi r_1 r_2}{r_1 + r_2} (\Bn \cdot \nabla H) (\Bp)
 \ee
and
 \be\label{estb2}
 \norm{\nabla b}_{L^{\infty} (\Om \setminus ( B_1 \cup B_2))} \leq C
 \ee
for a constant $C$ independent of $\ep$.
\end{thm}

It is worth looking more closely at the formula \eqnref{intensity2} of the stress intensity factor. The function $H$ is given by
 $$
 H(\Bx) = - \frac{1}{2\pi} \int_{\p \Om} \ln |\Bx-\By| g(\By) \, ds(\By) -
 \frac{1}{2\pi} \int_{\p \Om} \frac{\langle \Bx-\By, \nu(\By) \rangle}{|\Bx-\By|^2} \GL_\Om [g](\By) \, ds(\By),
 $$
and hence
 \begin{align}
 a &= -\frac {2\pi r_1 r_2}{\pi(r_1 + r_2)} \bigg[ \int_{\p \Om} \frac{\langle \Bp-\By, \Bn \rangle}{|\Bp-\By|^2} g(\By) \, ds(\By) \nonumber \\
 & \qquad + \int_{\p \Om} \Big(\frac{\langle \Bn, \nu(\By) \rangle}{|\Bp-\By|^2} - \frac{\langle \Bp-\By, \Bn \rangle \langle \Bp-\By, \nu(\By) \rangle}{|\Bp-\By|^4} \Big) \GL_\Om [g](\By) \, ds(\By) \bigg].
 \end{align}
So if we can measure the Dirichlet data $\GL_\Om [g]$ on $\p\Om$, we can determine the intensity of the stress using the boundary data. We emphasize that $a$ is bounded regardless of $\ep$ thanks to \eqref{GLg}.

%%%%%%%%%%%%%%%%%%%%%%%%%%%%%%%%%%%%%%%%%%%%%%%%%%%%%%%%%%%
\section{The insulated case}
%%%%%%%%%%%%%%%%%%%%%%%%%%%%%%%%%%%%%%%%%%%%%%%%%%%%%%%%%%%

We now deal with the case when circular inclusions are insulated, {\it i.e.}, the conductivities are $0$. Consider the solution to the free space problem:
\begin{equation} \label{eq:main_insul}
\begin{cases}
\ds\Delta u  = 0\quad&\mbox{in }{\mathbb{R}^2 \setminus \overline{B_1 \cup B_2}},\\
\ds \pd{u}{\nu}  = 0  \quad& \mbox{on }\p B_1\cup \p B_2,\\
\ds u(\Bx)- H(\Bx) = O(|\Bx|^{-1})\quad&\mbox{as } |\Bx| \rightarrow \infty.
\end{cases}
\end{equation}
From the jump formula of the single layer potential, $u$ can be represented as
 \be\label{ufreeinsul}
 u(\Bx) = H(\Bx) +\Scal_{B_1}[\varphi_1](\Bx) +\Scal_{B_2}[\varphi_2](\Bx),\quad \Bx\in \RR^2,
 \ee
for a pair of potentials $(\varphi_1,\varphi_2)\in L^2_0(\p B_1)\times L^2_0(\p B_2)$ satisfying \eqnref{varphi:eqn} with $\lambda = -\frac{1}{2}$.

Let $\widetilde{H}$ be an harmonic function in $\RR^2$ such that $H$ is a harmonic conjugate of $\widetilde{H}$. Then the solution $u$ to \eqnref{eq:main_insul} is a harmonic conjugate in $\RR^2\setminus\overline{B_1\cup B_2}$ of $\tilde u$ which is the solution to \eqnref{eq:main} with $\widetilde{H}$ in the place of $H$, see for example \cite{AKL}. Note that by the Cauchy-Riemann equation, the tangential derivative of $\tilde u$ is the same as the normal derivative of $u$ on the disks, and hence $\tilde u$ is constant on each disk $B_j$, $j=1.2$. Theorem \ref{Stiff} yields, for $\Bx \in \RR^2\setminus \overline{B_1 \cup B_2}$, as $\epsilon\rightarrow 0$,
 $$
 \nabla \tilde{u}(\Bx)= \frac {2 r_1 r_2}{r_1 + r_2} (\Bn \cdot \nabla \widetilde{H}) (\Bp)  \left( \frac {\mathbf{x}-\mathbf{p}_1} { \bigr|\mathbf{x}-\mathbf{p}_1\bigr|^2}- \frac {\mathbf{x}-\mathbf{p}_2} {\bigr|\mathbf{x}-\mathbf{p}_2\bigr|^2}\right) + O(1).
 $$
Let  $\Bt$ is the unit vector perpendicular to $\Bn$ such that $(\Bn, \Bt)$ is positively oriented and $\Bx^\perp = \begin{bmatrix} -x_2 \\ x_1 \end{bmatrix}$ for $\Bx\in\RR^2$.
Since $\nabla u = \left(\nabla \tilde u \right)^\perp$ and $\Bn \cdot \nabla \tilde{H}(\Bp) = \Bt \cdot \nabla H(\Bp)$, we have
 \be\label{uahb4}
 \nabla u(\Bx)= \frac {2 r_1 r_2}{r_1 + r_2} (\Bt \cdot \nabla H) (\Bp)  \left( \frac {(\mathbf{x}-\mathbf{p}_1)^\perp} { \bigr|\mathbf{x}-\mathbf{p}_1\bigr|^2}- \frac {(\mathbf{x}-\mathbf{p}_2)^\perp} {\bigr|\mathbf{x}-\mathbf{p}_2\bigr|^2}\right) + O(1).
 \ee

Using \eqnref{uahb4} we can obtain an expression of the solution $u$ to \eqnref{eq:main_insul}. Let $\arg:\mathbb{R}^2\setminus \{(0,0)\} \rightarrow [-\pi, \pi)$ be the argument function with a branch cut along the negative real axis, where $\Bx = (x_1,x_2)$ is identified with $x_1 + i x_2$. Define
 \be\label{hbot}
 h_\bot(\Bx) = \frac 1 {2\pi} \Bigr(\arg (\mathbf{x}-\mathbf{p}_1)-\arg (\mathbf{x}-\mathbf{p}_2) - \arg (\mathbf{x}-\mathbf{c}_1)+ \arg (\mathbf{x}-\mathbf{c}_2)\Bigr),
 \ee
where $\Bc_j$ is the center of $B_j$, $j=1,2$.  Note that $h_\bot$ is a harmonic function well defined in $\RR^2 \setminus \overline{(B_1\cup B_2)}$ since the jump discontinuity of the argument crossing the branch cut is canceled out owing to $\Bp_j, \Bc_j\in B_j$.
Moreover, we have
 \begin{align*}
 \nabla h_\bot(\Bx) &= \frac 1 {2\pi} \nabla \Bigr(\arg (\mathbf{x}-\mathbf{p}_1)-\arg (\mathbf{x}-\mathbf{p}_2) -  \arg (\mathbf{x}-\mathbf{c}_1)+ \arg (\mathbf{x}-\mathbf{c}_2) \Bigr) \\
 %&= \frac 1 {2\pi} (\nabla  h (\Bx))^\perp + O(1) \\
 &= \frac 1 {2\pi} \left( \frac {(\mathbf{x}-\mathbf{p}_1)^\perp} { \bigr|\mathbf{x}-\mathbf{p}_1\bigr|^2}- \frac {(\mathbf{x}-\mathbf{p}_2)^\perp} {\bigr|\mathbf{x}-\mathbf{p}_2\bigr|^2}\right) + O(1).
 \end{align*}

Similarly to the free space case, the solution $u$ to the boundary value problem with the insulated inclusion becomes the solution of the perfectly conducting disk by taking its conjugate. To be more precise, if $u$ is the solution
 \be\label{bvp-insul}
 \begin{cases}
 \ds  \Delta u=0 \quad &\mbox{in } \Om \setminus \overline{B_1 \cup B_2},\\
 \ds \pd{u}{\nu}  = 0  \quad& \mbox{on }\p B_1\cup \p B_2,\\
 \ds \pd{u}{\nu} = g&\mbox{on }\p \Om,
 \end{cases}
 \ee
where $\Om$ is a simply connected bounded domain $\mathcal{C}^2$-boundary and $g\in L^2_0(\p \Om)$,
 Then we have \eqnref{ufreeinsul} and \eqnref{varphi:eqn} with $\lambda=-\frac{1}{2}$ and
 \begin{equation*}
 H(\Bx) = -\Scal_\Om[g](\Bx) + \Dcal_\Om[u|_{\p\Om}](\Bx), \quad \Bx \in \Om.
 \end{equation*}
Since $\Om$ is a simply connected domain, the harmonic function $H$ admits a conjugate function $\widetilde{H}$ in $\Om$.
Similarly to in free space, there is a harmonic conjugate $\tilde u$ of $u$ in $\RR^2\setminus\overline{B_1\cup B_2}$ satisfies \eqnref{bvp} with a harmonic conjugate $\widetilde{H}$ in the place of $H$.

Thus, we have the following theorem.
\begin{thm} \label{Stiff-insul}
Let $u$ be either the solution to \eqnref{eq:main_insul} or the solution to \eqnref{bvp-insul} in which case $H$ is the function defined by \eqnref{defH}. Then, $u$ can be expressed as follows:
 \be\label{decom-insul}
 u(\Bx)= a_\bot h_\bot(\Bx) + b_\bot (\Bx), \quad \Bx \ \mbox{outside } B_1 \cup B_2
 \ee
where
 \be\label{abot}
 a_\bot = \frac {4\pi r_1 r_2}{r_1 + r_2} (\Bt \cdot \nabla H) (\Bp)
 \ee
and
 \be
 \norm{\nabla b_\bot}_{L^{\infty} (\Om \setminus ( B_1 \cup B_2))} \leq C
 \ee
for a constant $C$ independent of $\ep$.
\end{thm}

%%%%%%%%%%%%%%%%%%%%%%%%%%%%%%%%%%%%
\section{Numerical computations}
%%%%%%%%%%%%%%%%%%%%%%%%%%%%%%%%%%

In this section we compute numerically the solutions to \eqnref{eq:main} and \eqnref{eq:main_insul}. Computation of the solution in the presence of closely located inclusions with conductivity $k=0$ or $\infty$ is known to be a hard problem. To understand this difficulty, let us consider a standard way of computing the solution using the boundary integral method.

The solution to \eqnref{eq:main} can be represented as
 \be\label{uHS2}
 u(\Bx)= H(\Bx) + \Scal_{B_1}[\vp_1](\Bx) + \Scal_{B_2}[\vp_2](\Bx), \quad \Bx \in \RR^2 \setminus (B_1 \cup B_2),
 \ee
where $(\vp_1, \vp_2)$ is the solution to \eqnref{varphi:eqn} with $\Gl=\frac{1}{2}$.
We can compute $(\vp_1, \vp_2)$ numerically by discretizing \eqnref{varphi:eqn} with $M$ number of equi-spaced points on each disks $B_i$, $i=1,2$.  Let $\Bx_{i}^k$, $k=1,\dots,M$, be the nodal points on $\p B_i$ and set \be
 A= \begin{bmatrix}
        \Gl I_M & A_{12} \\
        A_{21} & \Gl I_M \\
      \end{bmatrix},
 \quad Y=\begin{bmatrix}
                      Y_1 \\
                      Y_2 \\
                \end{bmatrix},
 \ee
where
 \be\label{ystandard}
 Y_1 = \begin{bmatrix}
 -\pd{H}{\nu^{(1)}}(\Bx_1^1) \\
 \vdots \\
 -\pd{H}{\nu^{(1)}}(\Bx_1^M)
 \end{bmatrix}, \quad
 Y_2 = \begin{bmatrix}
 -\pd{H}{\nu^{(2)}}(\Bx_2^1) \\
 \vdots \\
 -\pd{H}{\nu^{(2)}}(\Bx_2^M)
 \end{bmatrix} ,
 \ee
and $A_{12}$ and $A_{21}$ are the evaluation of the kernel of $\pd{}{\nu^{(1)}}\Scal_{B_2}$ and $\pd{}{\nu^{(2)}}\Scal_{B_1}$ at nodes on $\p B_1$ and $\p B_2$, respectively.
We then obtain $(\varphi_1, \varphi_2)$ by solving
 \be \label{AXY}
 A \begin{bmatrix} \varphi_1 \\ \varphi_2 \end{bmatrix} = Y.
 \ee

\begin{figure}[h!]
\begin{center}
\epsfig{figure=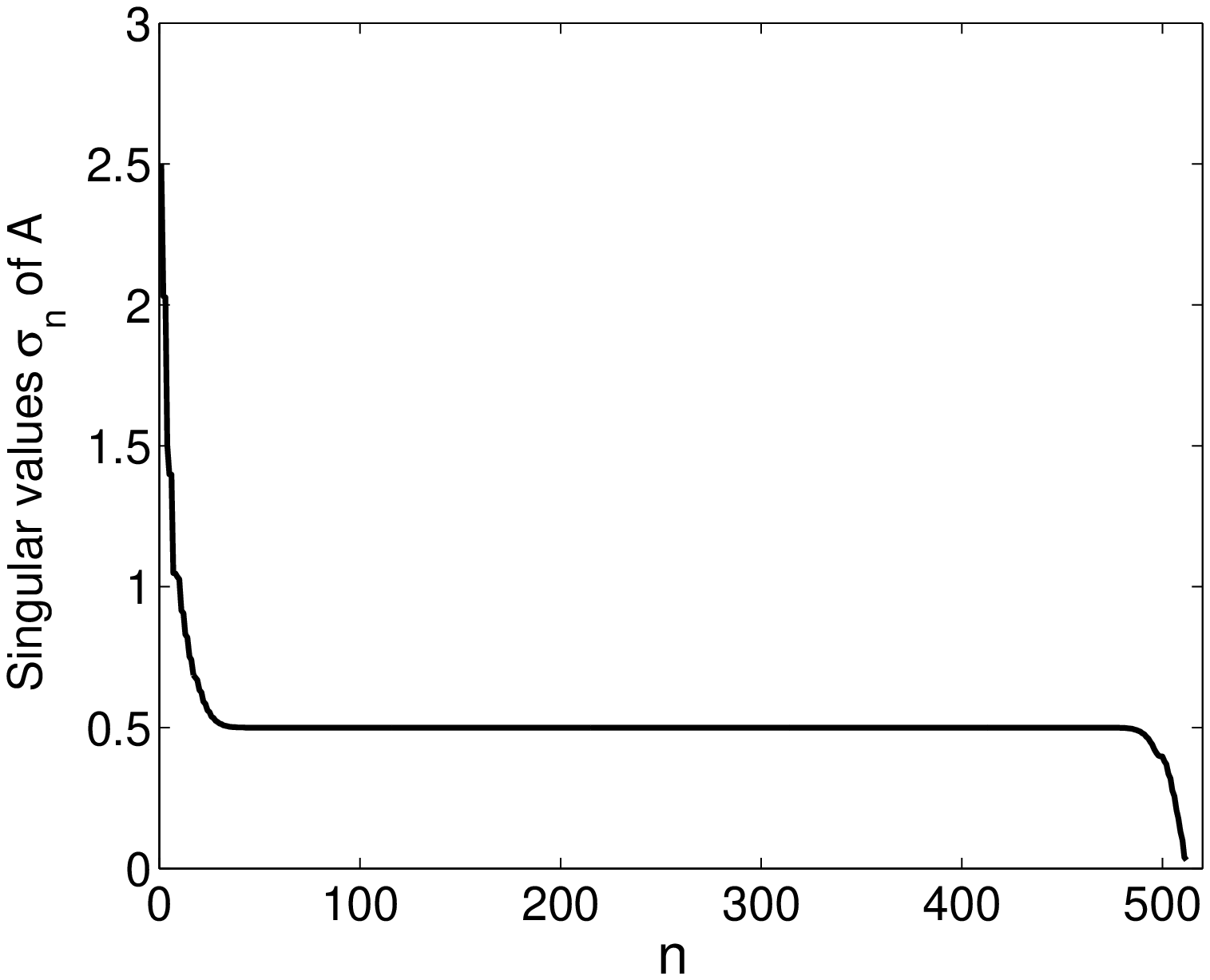,height = 5cm}\hskip .5cm
\epsfig{figure=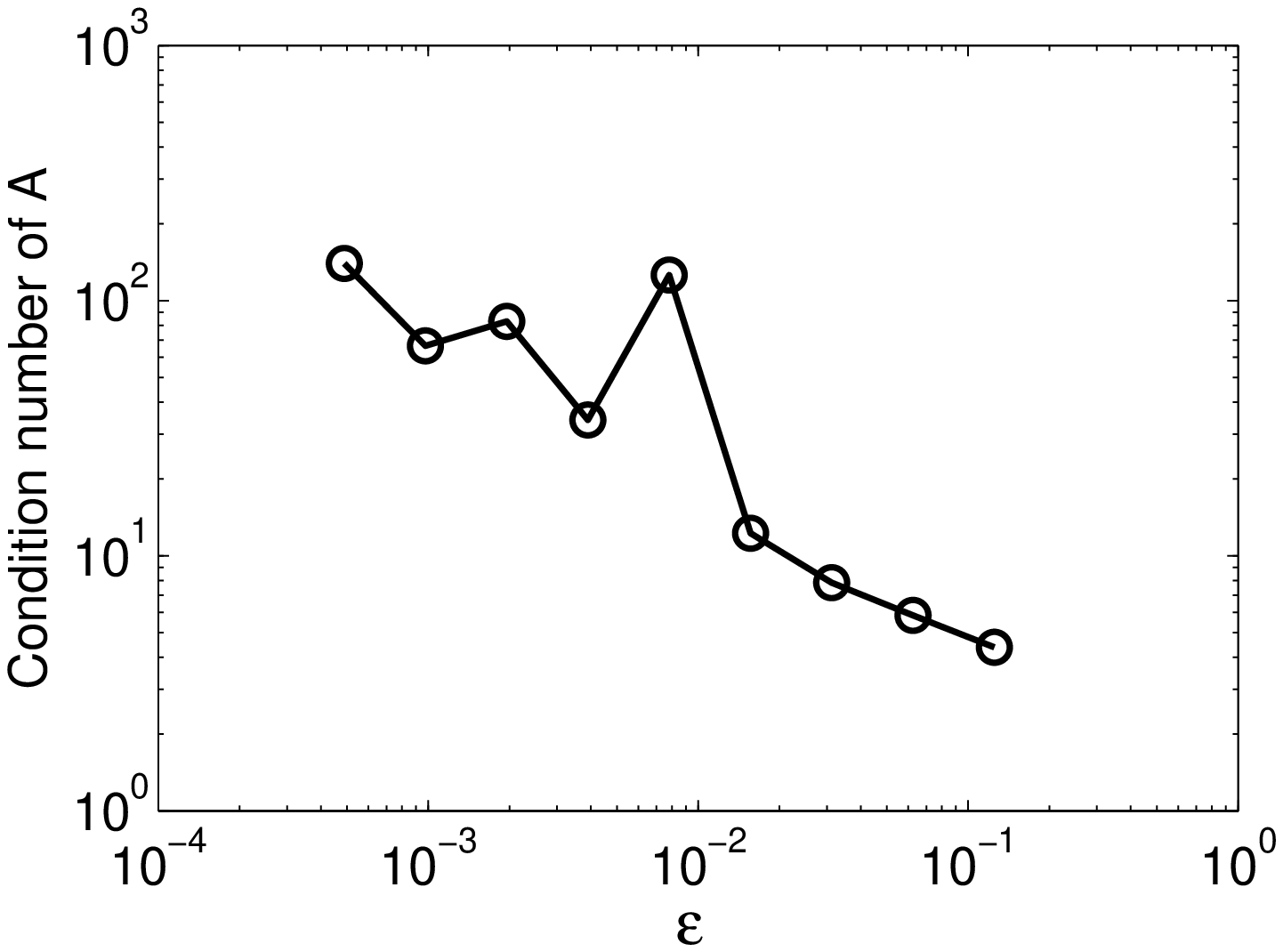, height = 5cm}
\end{center}
\caption{The first graph shows the singular values of $A$ in the decreasing order ($n$) when $\epsilon$ is 0.0020. The second graph shows the condition numbers of $A$ as $\ep$ tends to $0$ (from right to left).  We use 256 grid points on each $B_i$, $i=1,2$, and hence the dimension of $A$ is $512\times512$. }\label{singular_value}
\end{figure}

As Figure \ref{singular_value} shows, the matrix $A$ has small singular values, and the condition number of $A$ becomes worse as $\ep$ tends to $0$. Moreover, derivative of the kernel of $\pd{}{\nu}\Scal_{B_i}[\vp_i](\Bx)$, which is of the form $\frac{1}{2\pi}\frac{\la\Bx-\By,\nu(\Bx)\ra}{|\Bx-\By|^2}$, is as big as $\frac{1}{\epsilon^2}$ if $\Bx$ and $\By$ are on the arcs of $\p B_1$ and $\p B_2$ which are close to each other. Hence, if $\vp_i$ takes  large values on those arcs, the error in the discretization of the single layer potential becomes significant. From Theorem \ref{Stiff}, $\vp_i$ is as big as $\frac{1}{\sqrt{\epsilon}}$ when $\Bn \cdot \nabla H\neq 0$ at the middle point of the shortest line segment connecting $\p B_1$ and $\p B_2$. Therefore, we need finer grids as $\ep$ gets smaller, see Figure \ref{relative_error_fixedep}.

We will show that this difficulty can be overcome by using the characterization of singular terms given in \eqnref{uahb} and \eqnref{decom-insul}.

%%%%%%%%%%%%%%%%%%%%%%%%%%%%%%%%%%%%%%%%%%%%%%%%%%%%%%%%%%%%%
\subsection{Computation for the perfectly conducting case}
%%%%%%%%%%%%%%%%%%%%%%%%%%%%%%%%%%%%%%%%%%%%%%%%%%%%%%%%%%%%%

In this subsection we present a new method of computing the solution to \eqnref{eq:main} based on the characterization of the singular terms obtained in this paper.

Let
 \begin{align}
 \tilde{h}(\Bx) & = h(\Bx) - \frac{1}{2\pi}\bigr(\log |\Bx - \Bc_1| - \log |\Bx - \Bc_2| \bigr) \nonumber \\
 &= \frac{1}{2\pi}\bigr(\log |\Bx - \Bp_1| - \log |\Bx - \Bp_2| - \log |\Bx - \Bc_1| + \log |\Bx - \Bc_2| \bigr)
 \end{align}
for $\Bx\in\RR^2\setminus(B_1\cup B_2)$. This modified function has the property that $\int_{\p B_i} \tilde{h} =0$ for $i=1,2$, which is useful for the computation.

In view of \eqnref{uahb}, we look for a solution in the following form
 \be\label{perf-rep}
 u(\Bx) = a\tilde{h}(\Bx)+H(\Bx)+\Scal_{B_1}[\psi_1](\Bx)+\Scal_{B_2}[\psi_2](\Bx),\qquad \Bx\in\RR^2\setminus(B_1\cup B_2),
 \ee
instead of \eqnref{uHS2}, where $a$ is given by \eqnref{intensity}. According to Theorem \ref{Stiff}, the gradient of the function $H + \Scal_{B_1}[\psi_1] + \Scal_{B_2}[\psi_2]$ is bounded on $\Om\setminus(B_1\cup B_2)$ for some bounded set $\Om$ containing $B_1\cup B_2$, and hence $\| \psi_1 \|_{L^\infty(\p B_1)}$ and $\| \psi_2 \|_{L^\infty(\p B_2)}$ are bounded regardless of $\ep$.

To find the integral equation for density functions $(\psi_1, \psi_2)$, we argue as follows: Let $\tilde{h}^e$ be the extension of $\tilde{h}$ defined by $\tilde{h}^e (\Bx)=\tilde{h}(\Bx)$ for $\Bx \in \RR^2 \setminus (B_1 \cup B_2)$ and
 \be\label{tildeh_e}
  \tilde{h}^e(\Bx) =
  \begin{cases}
  \ds  h|_{\p B_1}-\frac{\log r_1}{2\pi} +\frac{1}{2\pi}\log |\Bx - \Bc_2|,\quad \Bx\in B_1, \\
  \nm
  \ds h|_{\p B_2}+\frac{\log r_2}{2\pi} -\frac{1}{2\pi}\log |\Bx - \Bc_1|,\quad \Bx\in B_2.
   \end{cases}
 \ee
Then $\tilde{h}^e=\tilde{h}$ on $\p B_i$ for $i=1,2$, and $\tilde{h}^e$ is harmonic in $B_1$ and $B_2$ as well as in $\RR^2 \setminus \overline{B_1 \cup B_2}$. Define
 $$
 u^e(\Bx) = a\tilde{h}^e(\Bx)+H(\Bx)+\Scal_{B_1}[\psi_1](\Bx)+\Scal_{B_2}[\psi_2](\Bx),\qquad \Bx\in\RR^2.
 $$
Then $u^e$ is continuous in $\RR^2$ and harmonic in $B_1\cup B_2$. Since $u^e = u$ is constant on $\p B_i$, $i=1,2$, $u^e$ is constant in $B_i$, $i=1,2$. By taking the interior normal derivative of $u^e$, one can see that $(\psi_1,\psi_2)\in L^2_0(\p B_1)\times L^2_0(\p B_2)$ is the unique solution to
 \be\label{psi_mod}
 \left\{ \begin{array}{l}
 \ds \frac{1}{2} \psi_1 +\pd{(\Scal_{B_2} [\psi_2] )}{\nu^{(1)}}
 = -\pd{H}{\nu^{(1)}}
- a \pd{\tilde{h}^e}{\nu^{(1)}}\bigr|_-
  \quad\mbox{on } \p B_1, \\
 \nm
 \ds \pd{(\Scal_{B_1} [\psi_1] )}{\nu^{(2)}} + \frac{1}{2} \psi_2 = -\pd{H}{\nu^{(2)}}
   - a \pd{\tilde{h}^e}{\nu^{(2)}}\bigr|_-
 \quad\mbox{on } \p B_2.
 \end{array}  \right.
 \ee
Moreover, we have from \eqnref{tildeh_e}
\begin{align*}
\pd{\tilde{h}^e}{\nu^{(1)}}\bigr|_- (\Bx) &=
\frac{1}{2\pi}\frac{\la \Bx - \Bc_2,\nu^{(1)}(\Bx)\ra}{|\Bx - \Bc_2|^2}, \quad \Bx \in \p B_1, \\
\pd{\tilde{h}^e}{\nu^{(2)}}\bigr|_- (\Bx) &= - \frac{1}{2\pi}\frac{\la \Bx - \Bc_1,\nu^{(2)}(\Bx)\ra}{|\Bx - \Bc_1|^2},  \quad \Bx \in \p B_2.
\end{align*}

We can discretize \eqnref{psi_mod} and solve \eqnref{AXY} to obtain $(\psi_1,\psi_2)$. Here $Y_1$ and $Y_2$ are given by
 \be\label{ynew}
 Y_1 = \begin{bmatrix} -\pd{H}{\nu^{(1)}}(\Bx_1^1)- \frac{a}{2\pi}\frac{\la \Bx_1^1 - \Bc_2,\ \nu^{(1)}(\Bx_1^1)\ra}{|\Bx_1^1 - \Bc_2|^2} \\
 \vdots \\
 -\pd{H}{\nu^{(1)}}(\Bx_1^M)- \frac{a}{2\pi}\frac{\la \Bx_1^M - \Bc_2,\ \nu^{(1)}(\Bx_1^M)\ra}{|\Bx_1^M - \Bc_2|^2}
 \end{bmatrix}
, \quad
 Y_2 = \begin{bmatrix} -\pd{H}{\nu^{(2)}}(\Bx_2^1)+ \frac{a}{2\pi}\frac{\la \Bx_2^1 - \Bc_1,\ \nu^{(2)}(\Bx_2^1)\ra}{|\Bx_2^1 - \Bc_1|^2} \\
 \vdots \\
 -\pd{H}{\nu^{(2)}}(\Bx_2^M)+ \frac{a}{2\pi}\frac{\la \Bx_2^M - \Bc_1,\ \nu^{(2)}(\Bx_2^M)\ra}{|\Bx_2^M - \Bc_1|^2}
 \end{bmatrix} .
 \ee

While $(\vp_1, \vp_2)$ in the representation \eqnref{uHS2} increases arbitrarily as $\ep$ tends to 0, $(\psi_1,\psi_2)$ stays bounded. The difference between the actual $(\psi_1, \psi_2)$ and the numerically obtained one is much smaller than that for $(\vp_1,\vp_2)$ as Figure \ref{relative_error} shows.

\begin{figure}[h!]
\begin{center}
\epsfig{figure=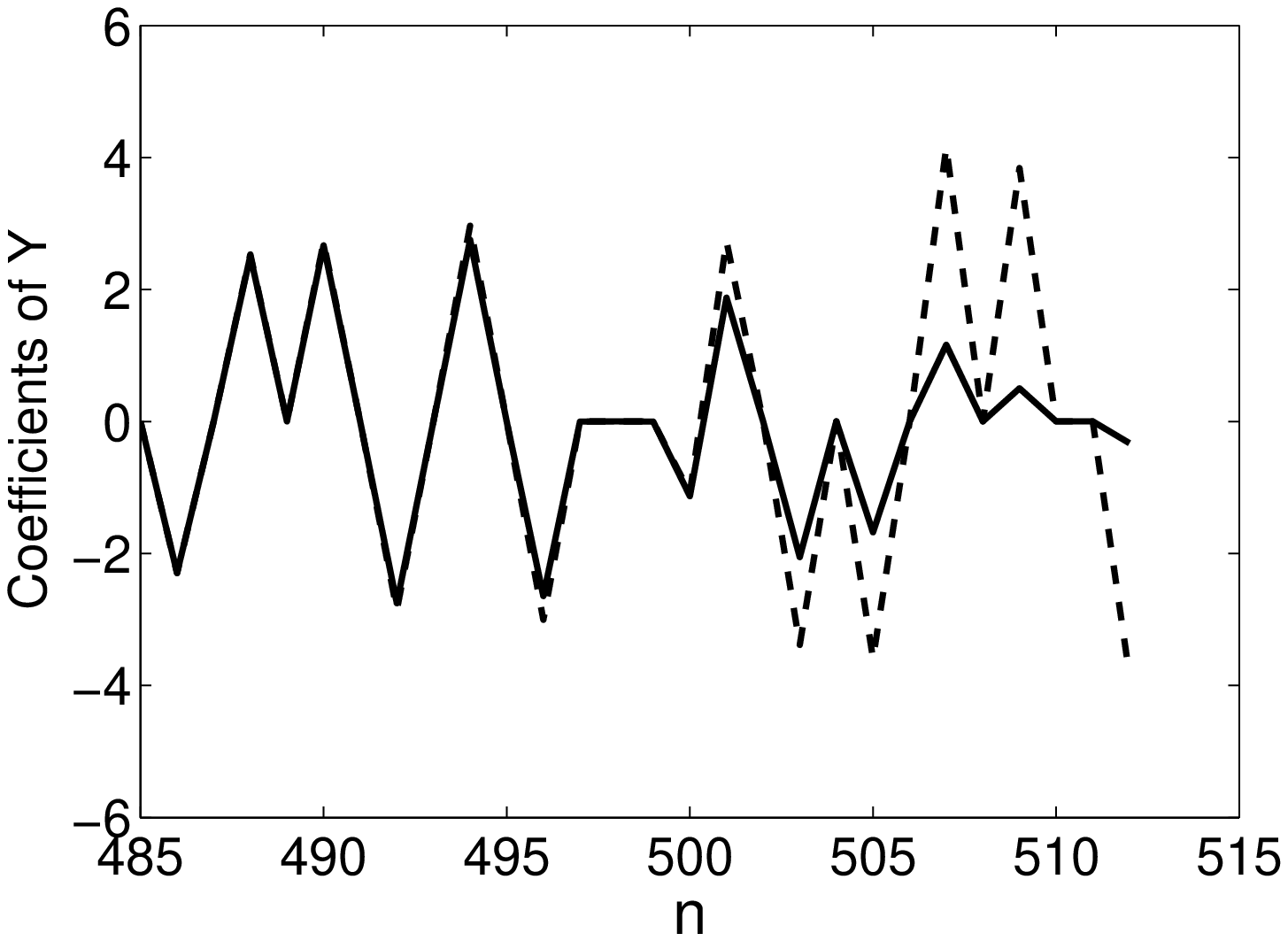, height = 5cm} \hskip .5cm
\epsfig{figure=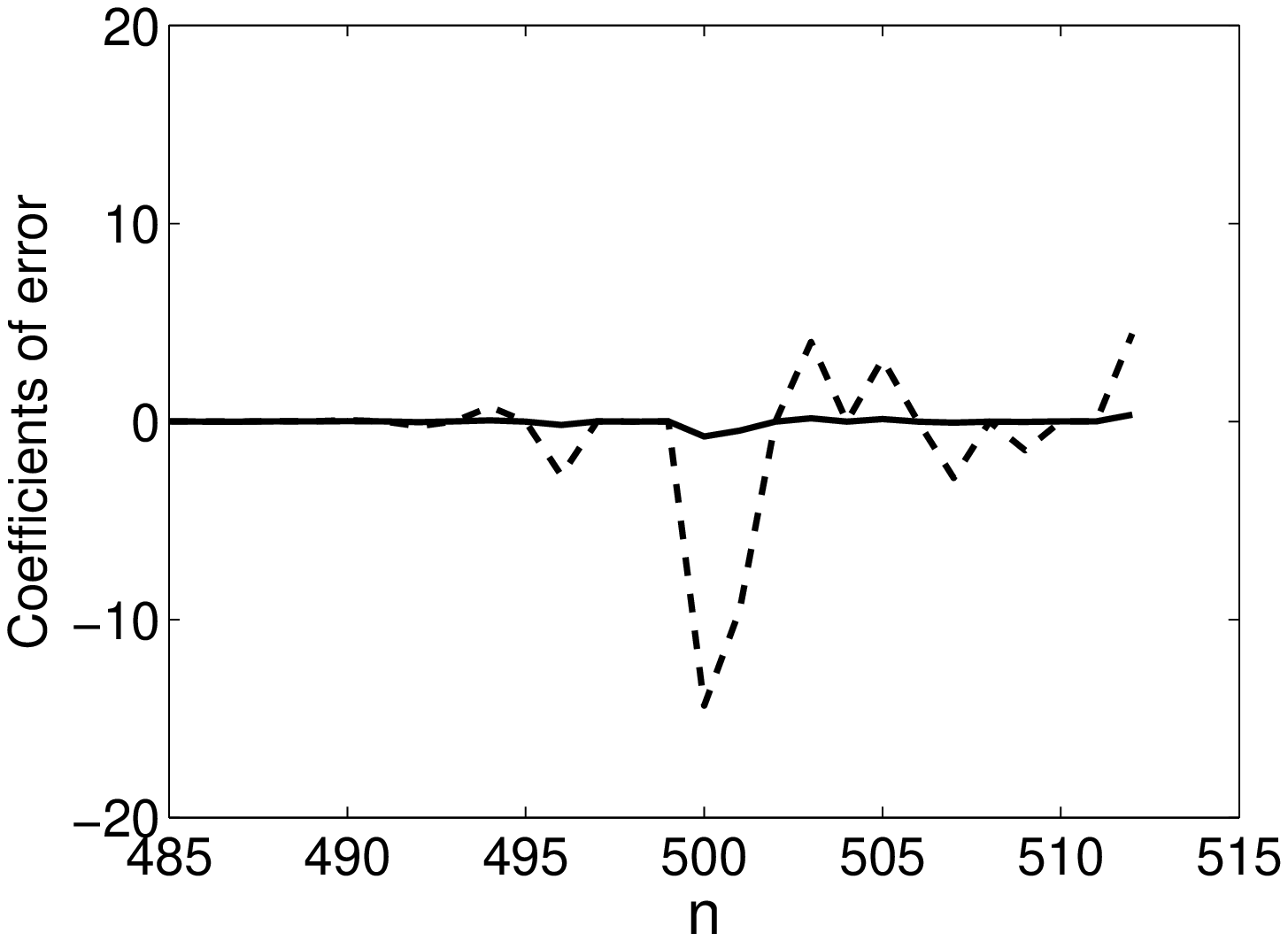, height = 5cm}
\end{center}
\caption{The first graph shows the inner product of singular vectors (corresponding small singular values) of $A$ with $Y$ in \eqnref{AXY}, where $Y$ is the discretization of the right hand side of \eqnref{varphi:eqn} (the dashed line), and of \eqnref{psi_mod} (the solid line). The second graph shows the inner product of singular vectors (corresponding small singular values) of $A$ with the error in \eqnref{AXY}. The radii of disks are fixed as $r_1= r_2=1$, and the number of equi-spaced grid points is 256 on each $\p B_i$. The distance $\ep=0.0020$, and the background potential is given by $H(\Bx) =x_1$. $n$ indicates the location of singular values when listed in decreasing order.}\label{coefficient2}
\end{figure}

The first graph of Figure \ref{coefficient2} shows the inner products of $Y$ in \eqnref{AXY} with singular vectors of $A$ corresponding to small singular values. The dotted graph is when \eqnref{ystandard} is used and solid one is when \eqnref{ynew} is used. The inner product using \eqnref{ystandard} is larger than the one using \eqnref{ynew}. This is expected: since the difference between the single layer potential and its discretization is large in the narrow region between two disks, the  singular vector (corresponding to small singular values) components of this difference is not small. The second graph in Figure \ref{coefficient2} shows the inner products of $A \begin{bmatrix} \varphi_1 \\  \varphi_2 \end{bmatrix} -Y$ with singular vectors of $A$ corresponding to small singular values, where $(\vp,\vp_2)$ is the (numerical) solution to \eqnref{AXY} using $Y$ in \eqnref{ystandard} (dotted line) and in \eqnref{ynew} (solid line). This numerical solution is obtained using the method described in subsection \ref{5.3} (with high precision). The graphs clearly show that the method of this paper works much better than the standard method. Here  $H(\Bx)=x_1$.

%%%%%%%%%%%%%%%%%%%%%%%%%%%%%%%%%%%%%%%%%%%%%%
\subsection{Computation for the insulated case}
%%%%%%%%%%%%%%%%%%%%%%%%%%%%%%%%%%%%%%%%%%%%%%

Let $h_\bot$ be the function defined by \eqnref{hbot}.
Since $\arg(\Bx - \Bp_1) - \arg(\Bx - \Bp_2)- \arg(\Bx - \Bc_1)$ is a harmonic conjugate of $\log|\Bx-\Bp_1|-\log|\Bx-\Bp_2|-\log|\Bx-\Bc_1|$, which is constant on $\p B_1$, we have
 $$
 \pd{}{\nu}\bigr(\arg(\Bx - \Bp_1) - \arg(\Bx - \Bp_2)- \arg(\Bx - \Bc_1) \bigr)=0,\quad \mbox{on }\p B_1.
 $$
Hence,
\be\label{hcsys}
\pd{h_\bot}{\nu^{(1)}}(\Bx)=\frac{1}{2\pi}\pd{(\arg(\Bx -\Bc_2))}{\nu^{(1)}}\Bigr|_{\p B_1} = \frac{1}{2\pi}\frac{\la (\Bx - \Bc_2)^\perp,\nu^{(1)}(\Bx)\ra}{|\Bx - \Bc_2|^2}, \quad \Bx \in \p B_1.
\ee
Similarly, we have
\be
\ds\pd{h_\bot}{\nu^{(2)}} (\Bx) =-\frac{1}{2\pi}\pd{(\arg(\Bx - \Bc_1))}{\nu^{(2)}}\Bigr|_{\p B_2}=- \frac{1}{2\pi}\frac{\la (\Bx - \Bc_1)^\perp,\nu^{(2)}(\Bx)\ra}{|\Bx - \Bc_1|^2}, \quad \Bx \in \p B_2.
\ee

We look for a solution $u$ to \eqnref{eq:main_insul} in the following form:
 \be\label{insul-rep}
 u(\Bx) = a_\bot h_\bot(\Bx)+H(\Bx)+\Scal_{B_1}[\psi_1](\Bx)+\Scal_{B_2}[\psi_2](\Bx),\qquad x\in\RR^2\setminus(B_1\cup B_2),
 \ee
where $a_\bot$ is given by \eqnref{abot} and $(\psi_1,\psi_2)\in L^2_0(\p B_1)\times L^2_0(\p B_2)$ is the solution to
 \be\label{psi_insu_mod}
 \left\{ \begin{array}{l}
 \ds \frac{1}{2} \psi_1 -\pd{(\Scal_{B_2} [\psi_2] )}{\nu^{(1)}}
 = \pd{H}{\nu^{(1)}}
 + \frac{a}{2\pi}\frac{\la (\Bx - \Bc_2)^\perp,\nu^{(1)}(\Bx)\ra}{|\Bx - \Bc_2|^2}
  \quad\mbox{on } \p B_1, \\
 \nm
 \ds -\pd{(\Scal_{B_1} [\psi_1] )}{\nu^{(2)}} + \frac{1}{2} \psi_2 = \pd{H}{\nu^{(2)}}
  - \frac{a}{2\pi}\frac{\la (\Bx - \Bc_1)^\perp,\nu^{(2)}(\Bx)\ra}{|\Bx - \Bc_1|^2}
 \quad\mbox{on } \p B_2.
 \end{array}  \right.
 \ee

%%%%%%%%%%%%%%%%%%%%%%%%%%%%%%%%%%%%
\subsection{Numerical Illustration}\label{5.3}
%%%%%%%%%%%%%%%%%%%%%%%%%%%%%%%%%%%%
In this subsection, we illustrate results of numerical computations using the algorithms proposed in the previous subsections. Two discs are $B_j = B(\Bc_j, r_j)$, $j=1,2$, of radius $r_j$ and centered at $\Bc_1=(-r_1-\epsilon/2, 0)$ and $\Bc_2 = (r_2+\epsilon/2, 0)$.

We compute the solution in two different ways and compare them to demonstrate the effectiveness of the method proposed in this paper. We first compute the solution using the standard representation of the solution, namely, we use \eqnref{uHS2} and solve numerically \eqnref{varphi:eqn}. The discretization for the computation was described at the beginning of this section. We denote by $u$ the solution computed by this method. We then compute the solution using the representation \eqnref{perf-rep} and solve \eqnref{psi_mod}. The solution is denoted by $u^h$. For comparison we solve \eqnref{psi_mod} yet another method which provides the solution with higher precision (but with high cost).

\begin{figure}[h!]
\begin{center}
\epsfig{figure=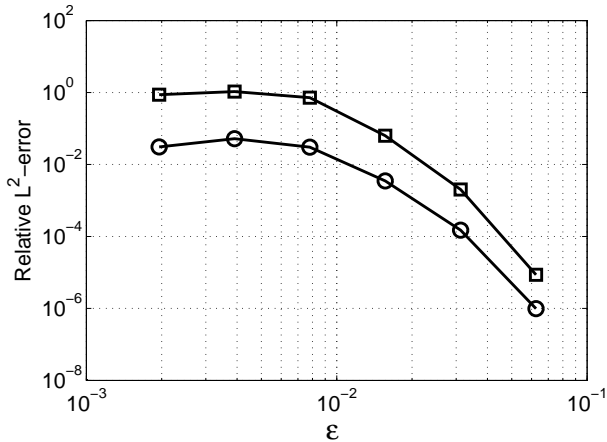,height = 6cm}
\end{center}
\caption{The circles are the relative errors of $\pd{u^h}{\nu}$ compared to $\pd{u^R}{\nu}$, and the squares are those  of $\pd{u}{\nu}$. Both depend on the distance $\epsilon$ between two perfectly conducting disks. The solution based on the asymptotic expansion has much smaller error. We use 256 grid points on each $B_i$, $i=1,2$. }\label{relative_error}
\end{figure}

Let $R_i$, $i=1,2$, be the reflection with respect to the disks $B_i$ defined by \eqnref{ref_B}.
We also define the the reflection of a function $f$ by
$(R_if)(\Bx) = f(R_i(\Bx))$ for $\Bx\in\RR^2$. Using the same argument as in \cite{AKL} (see also \cite{CG}), one can show that the solution $(\psi_1,\psi_2)$ to \eqnref{psi_mod} is given by
\begin{align*}
\psi_1 = -2\sum_{m=0}^\infty\pd{}{\nu^{(1)}}\bigg[(R_2R_1)^m\Bigr(H+\frac{a}{2\pi}\log|\Bx-\Bc_2|
-R_2\bigr[H-\frac{a}{2\pi}\log|\Bx-\Bc_1|\bigr]\Bigr)\bigg]\Bigr|_{\p B_1}, \\
\psi_2 = -2\sum_{m=0}^\infty\pd{}{\nu^{(2)}}\bigg[(R_1R_2)^m(H-\frac{a}{2\pi}\log|\Bx-\Bc_1|
-R_1\bigr[H+\frac{a}{2\pi}\log|\Bx-\Bc_2|\bigr]\Bigr)\bigg]\Bigr|_{\p B_2}.
\end{align*}
%\end{lem}
%\pf
%For a function $v$ harmonic in $\overline{B_i}$, $i=1,2$, we have
%\be
%\Scal_{B_i}\left[\pd{v}{\nu^{(i)}}\right](\Bx)=\frac{1}{2}R_iv(\Bx)+\frac{1}{2}v(\Bc_i),\quad\Bx\in\RR^2\setminus\overline{B_i}.
%\ee
%See
%Note that $(R_1R_2)^m(H-\frac{a}{2\pi}\log|\Bx-\Bc_1|-R_1\bigr[H+\frac{a}{2\pi}\log|\Bx-\Bc_2|\bigr ]\Bigr))$ is %harmonic in $\overline{B_2}$. We have
%\begin{align*}
%&\pd{(\Scal_{B_2}[\psi_2])}{\nu^{(1)}}\Bigr|_{\p B_1}\\
%&=-2\sum_{m=0}^\infty\pd{}{\nu^{(1)}}\Scal_{B_2}\Bigg[\pd{}{\nu^{(2)}}\bigg[(R_1R_2)^m(H-\frac{a}{2\pi}\log|\Bx-\Bc_1|
%-R_1\bigr[H+\frac{a}{2\pi}\log|\Bx-\Bc_2|\bigr]\Bigr)\bigg]\Bigr|_{\p B_2}\Bigg]\Bigr|_{\p B_1}\\
%&=-\sum_{m=0}^\infty\pd{}{\nu^{(1)}}\bigg[R_2(R_1R_2)^m(H-\frac{a}{2\pi}\log|\Bx-\Bc_1|
%-R_1\bigr[H+\frac{a}{2\pi}\log|\Bx-\Bc_2|\bigr]\Bigr)\bigg]\Bigr|_{\p B_1}\\
%&=\psi_1-\pd{}{\nu^{(1)}}\Bigr(H+\frac{a}{2\pi}\log|\Bx-\Bc_2|\Bigr).
%\end{align*}
%\qed
We denote by $u^R$ the solution obtained by this method. We compare these solutions for various values of $\epsilon$. The radii are fixed as $r_1= 2$ and $r_2 = 1.5$, and the number of equi-spaced grid points is 256 on each $\p B_i$. The background potential is given by $H(\Bx) = 2x_1 + (x_1^2 - x_2^2)$.
Figure \ref{relative_error} shows the relative $L^2$-errors of the normal flux $\pd{u}{\nu }$ and $\pd{u^h}{\nu }$, compared to the normal flux $\pd{u^R}{\nu }$. The vertical axis in the figure represents the values of
 $$
 \frac{ \left\| \pd{u^h}{\nu^{(1)}} - \pd{u^R}{\nu^{(1)}} \right\|_{L^2(\p B_1)}} {2 \left\| \pd{u^R}{\nu^{(1)}} \right\|_{L^2(\p B_1)}} + \frac{ \left\| \pd{u^h}{\nu^{(2)}} - \pd{u^R}{\nu^{(2)}} \right\|_{L^2(\p B_2)}} {2 \left\| \pd{u^R}{\nu^{(2)}} \right\|_{L^2(\p B_2)}}
 $$
indicated by circles and
 $$
 \frac{ \left\| \pd{u}{\nu^{(1)}} - \pd{u^R}{\nu^{(1)}} \right\|_{L^2(\p B_1)}} {2 \left\| \pd{u^R}{\nu^{(1)}} \right\|_{L^2(\p B_1)}} + \frac{ \left\| \pd{u}{\nu^{(2)}} - \pd{u^R}{\nu^{(2)}} \right\|_{L^2(\p B_2)}} {2 \left\| \pd{u^R}{\nu^{(2)}} \right\|_{L^2(\p B_2)}}
 $$
indicated by squares. As $\epsilon$ decreases (from right to left), the relative error increases for $u$ and $u^h$ as expected. However, the relative error of $u^h$ is notably small compared to that of $u$: when $\epsilon\sim 0.001$, the relative error of $u^h$ is as small as 0.01, but that of $u$ is as big as 1.

\begin{figure}[h!]
\begin{center}
\epsfig{figure=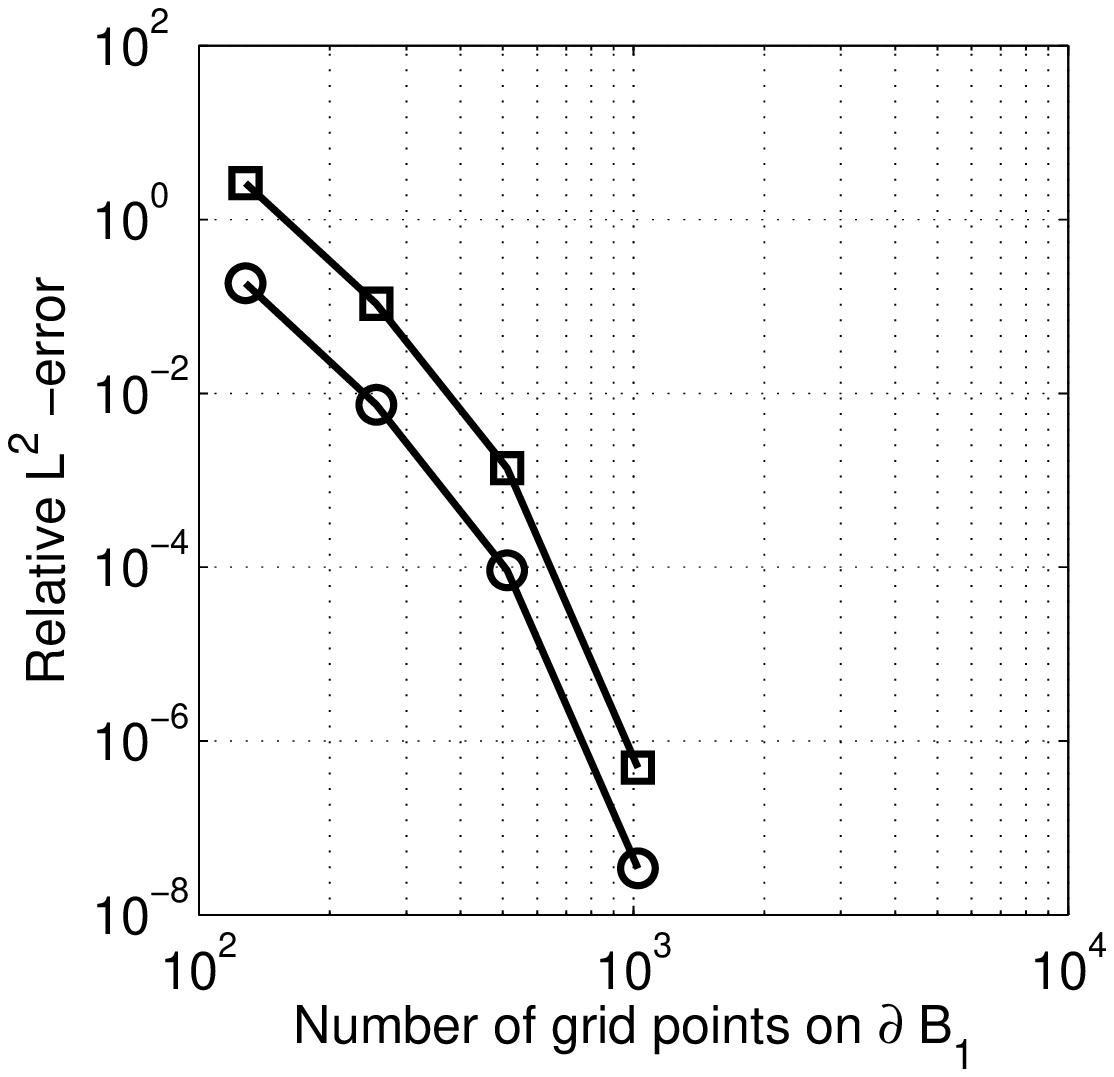, height=5cm}\hskip 1cm
\epsfig{figure=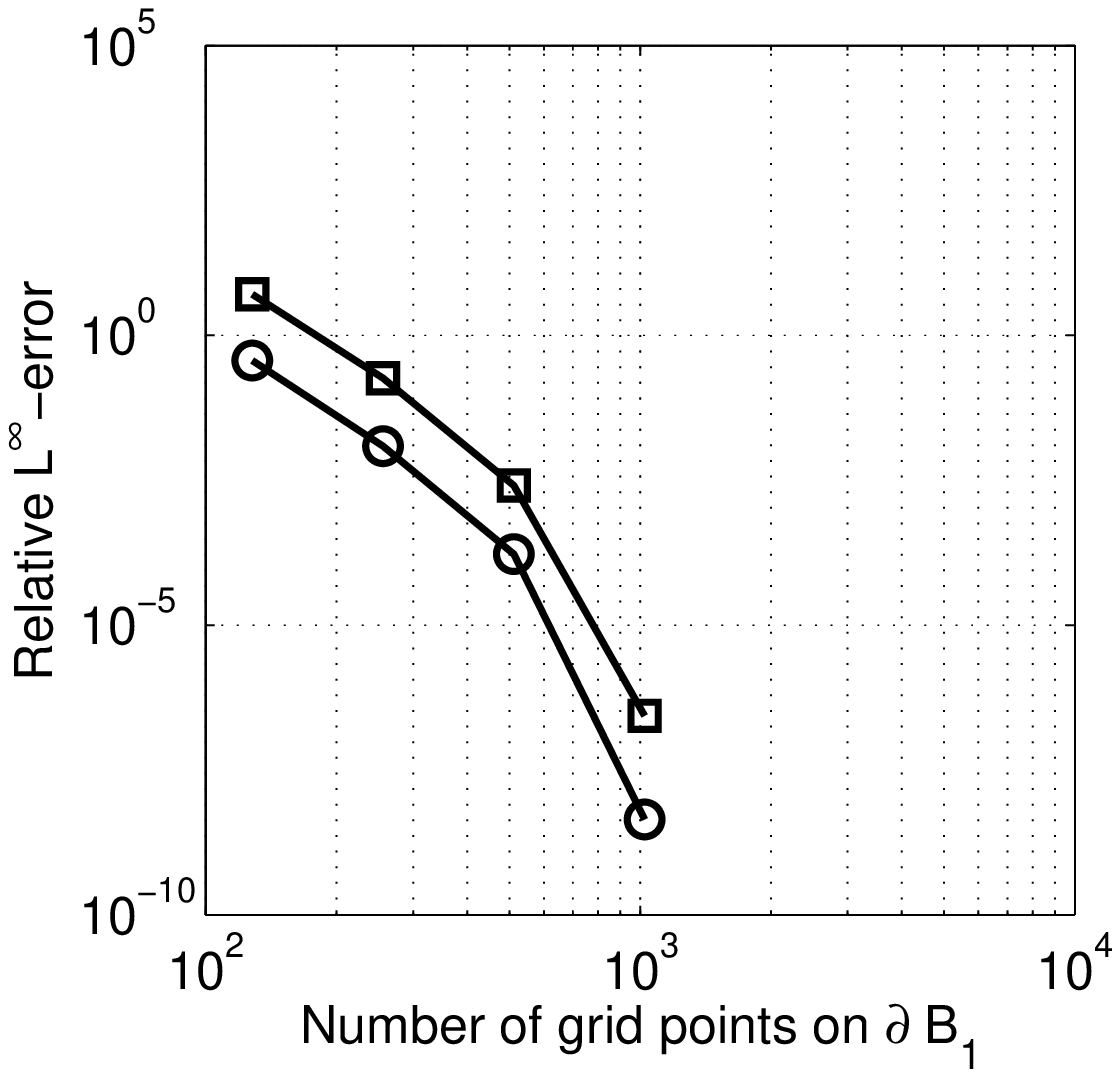,height=5cm}\\[.5cm]
\end{center}
\caption{The circles are the relative errors of $\pd{u^h}{\nu}$ compared to $\pd{u^R}{\nu}$, and the squares are those of $\pd{u^h}{\nu}$ in $L^2$ and $L^\infty$ norms for various grid numbers. The solution based on the asymptotic expansion has much smaller error.}\label{relative_error_fixedep}
\end{figure}

In Figure \ref{relative_error_fixedep}, fixing $\epsilon= 0.0156$, we compare the relative errors for different grid numbers.
Both radii are 1, and $H(\Bx) = x_1$. The difference between the normal flux $\pd{u}{\nu^{(j)} }$ and $\pd{u^R}{\nu^{(j)} }$ takes its maximum value at the point nearest the middle point of the shortest line segment between $\p B_1$ and $\p B_2$, which is the same for  $\pd{u^h}{\nu^{(j)} }$. As the grid numbers increase, both the relative $L^2$-error and the maximal difference decrease. But, relative errors of $u^h$ are much smaller than those of $u$.

\begin{figure}[h!]
\begin{center}
\epsfig{figure=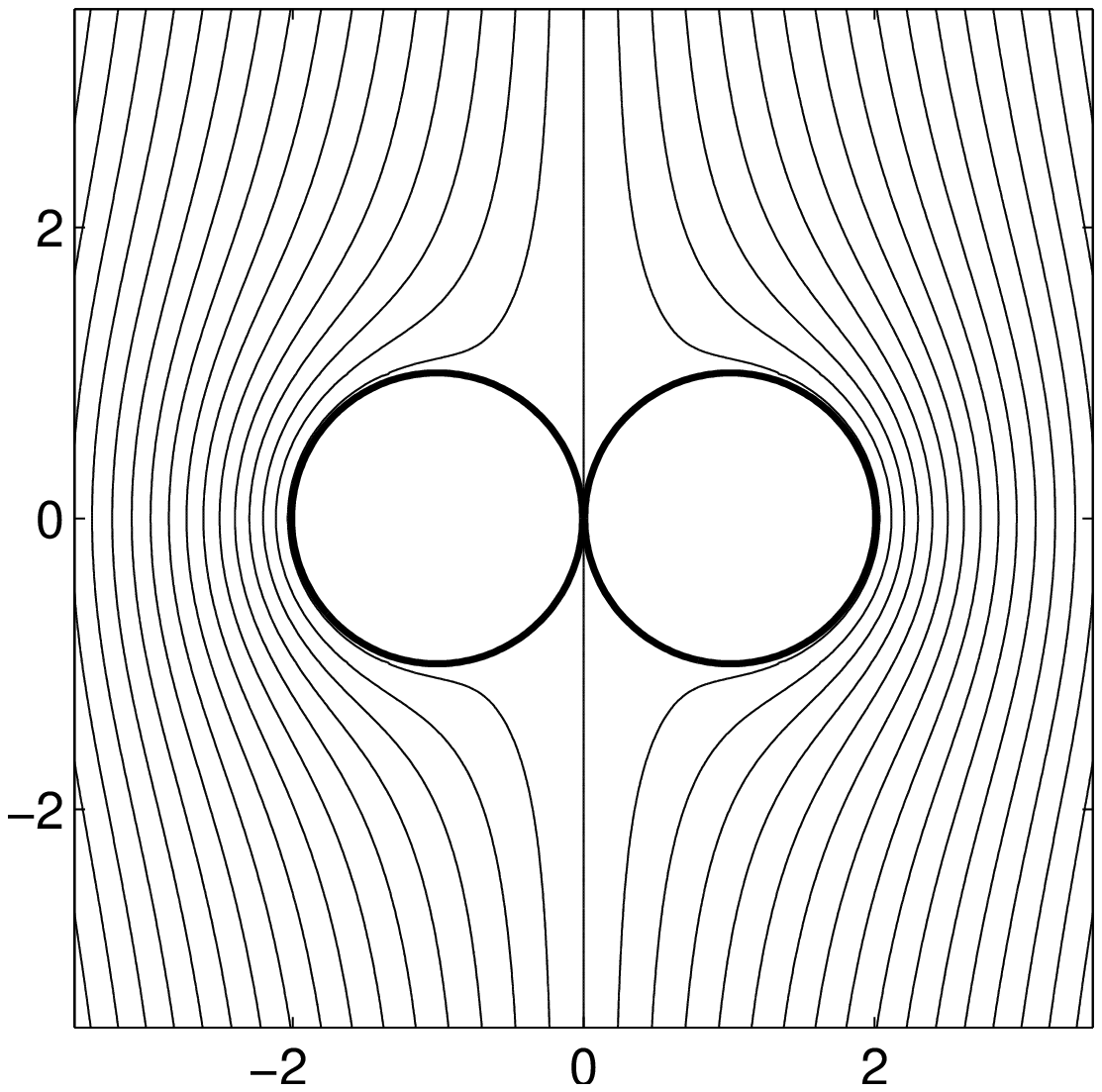, height=4cm}\hskip 1cm
\epsfig{figure=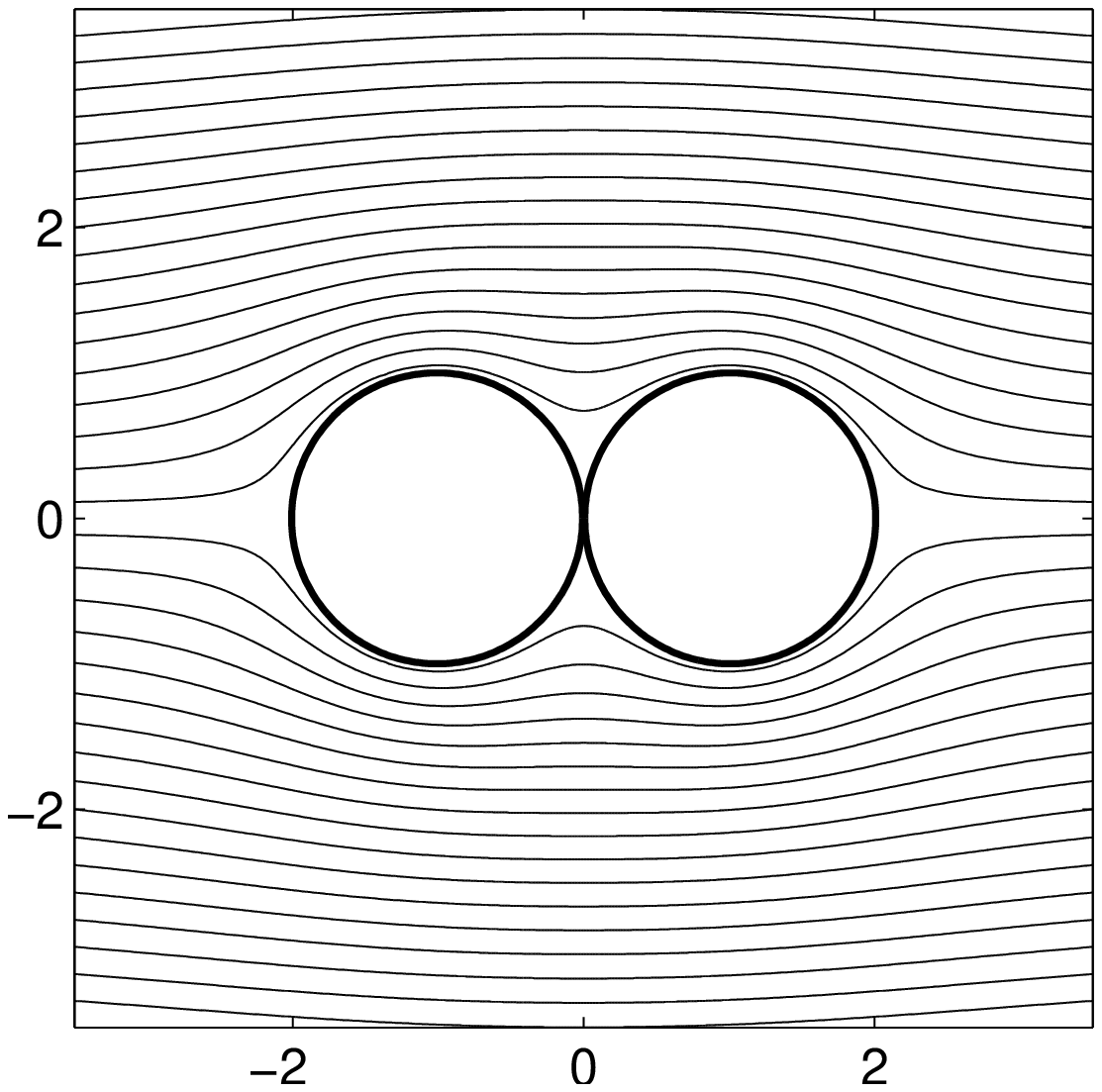, height=4cm}\\[.5cm]
\epsfig{figure=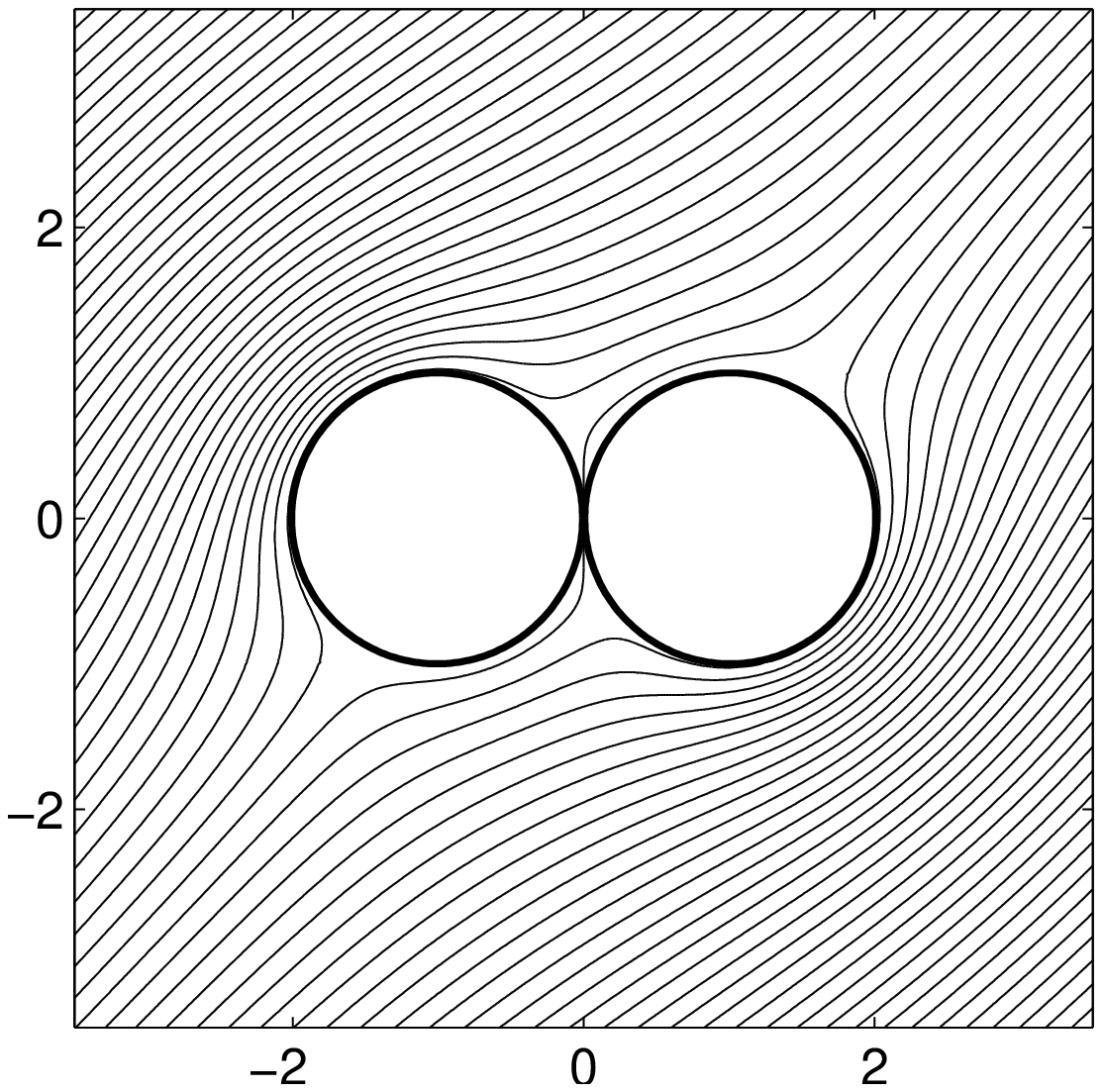, height=4cm}\hskip 1cm
\epsfig{figure=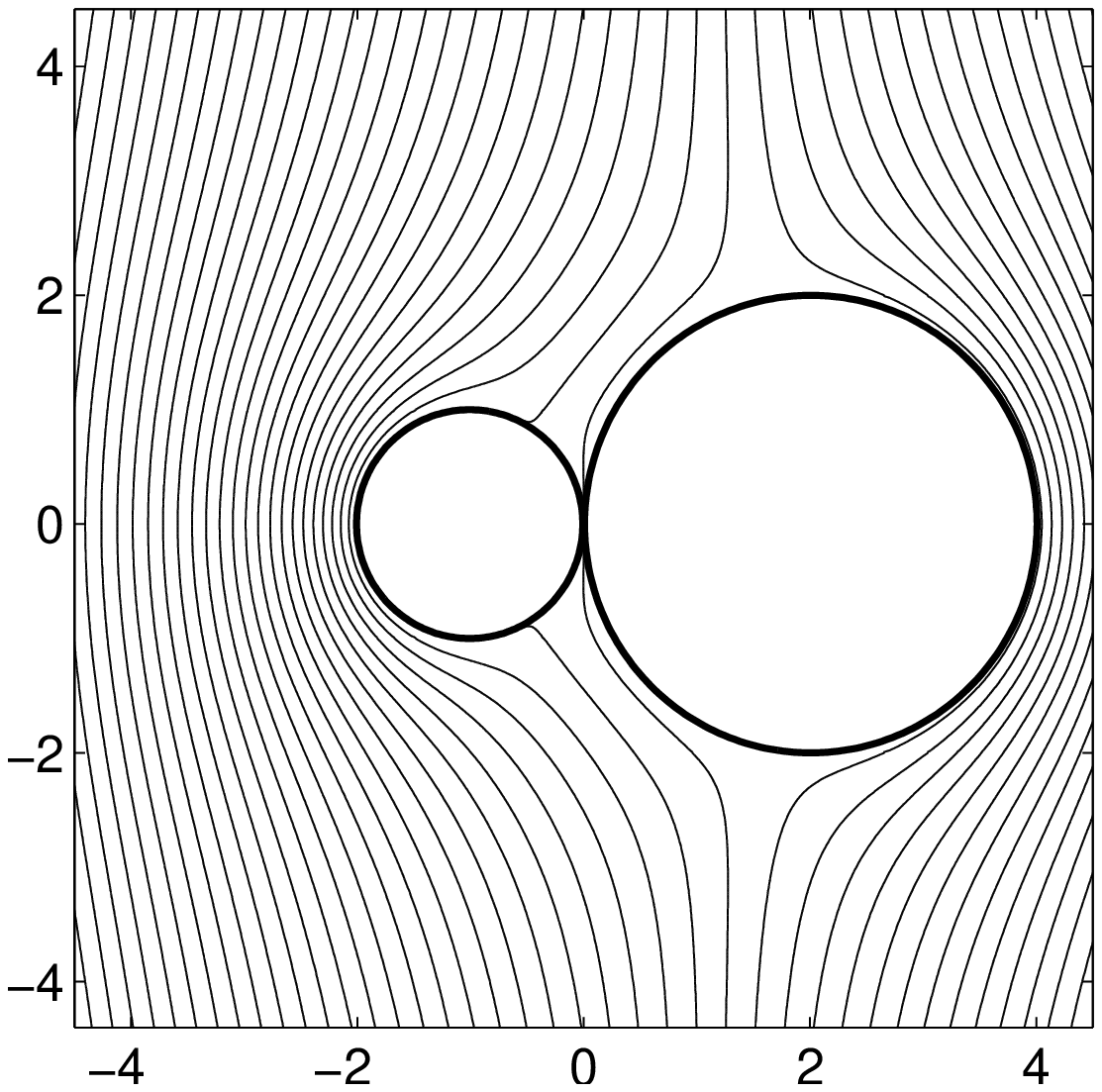, height=4cm}\\[.5cm]
\end{center}
\caption{Level curves of the free space conducting case. The entire harmonic function $H(\Bx) = x_1$ in the upper-left and the lower-right figure, and $H(\Bx) = x_2$ in the upper-right one, and
$H(\Bx) =x_1 -x_2$ in the lower-left one.}\label{level_curve}
\end{figure}

\begin{figure}[h!]
\begin{center}
\epsfig{figure=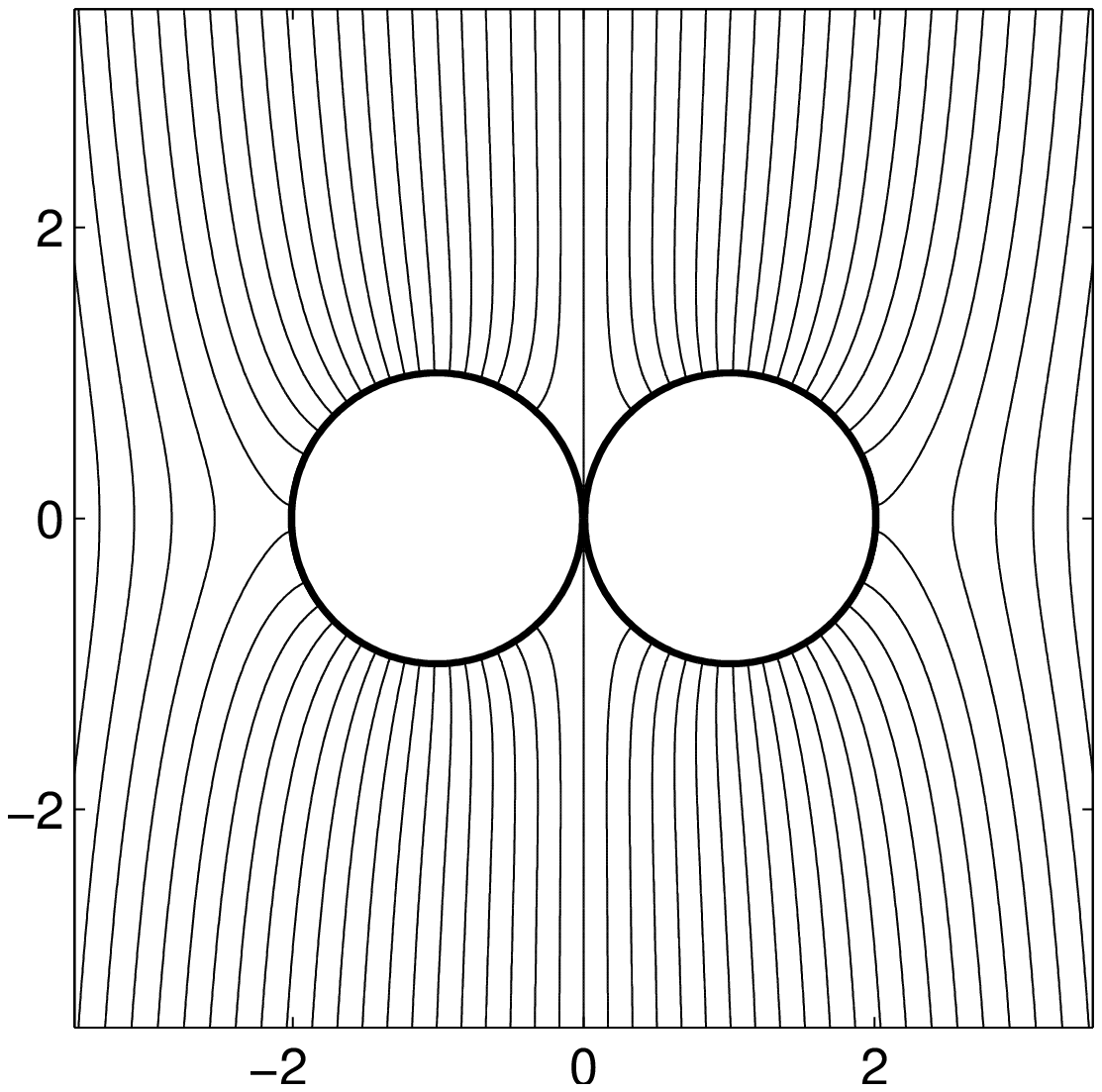, height=4cm}\hskip 1cm
\epsfig{figure=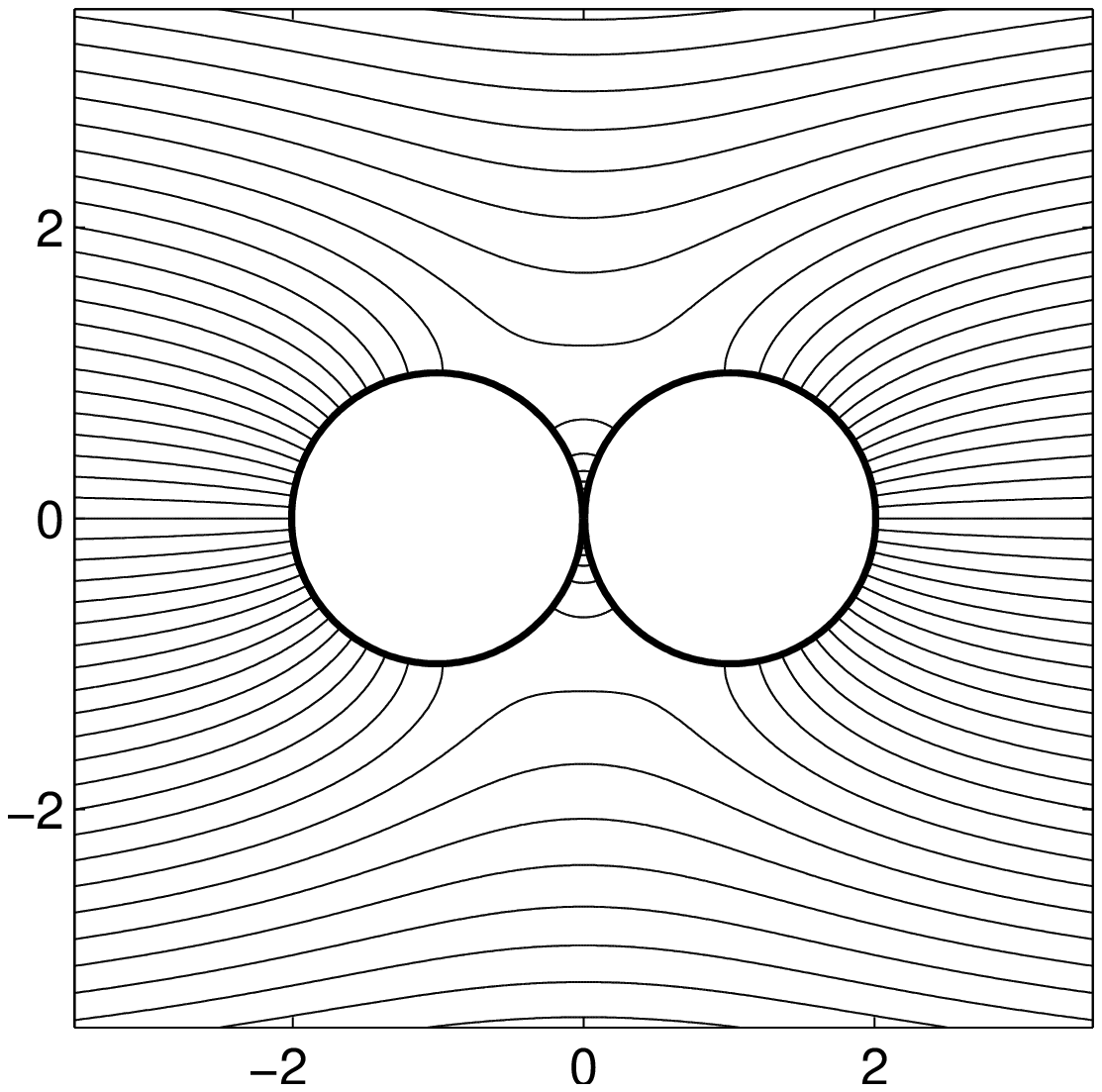, height=4cm}\\[.5cm]
\epsfig{figure=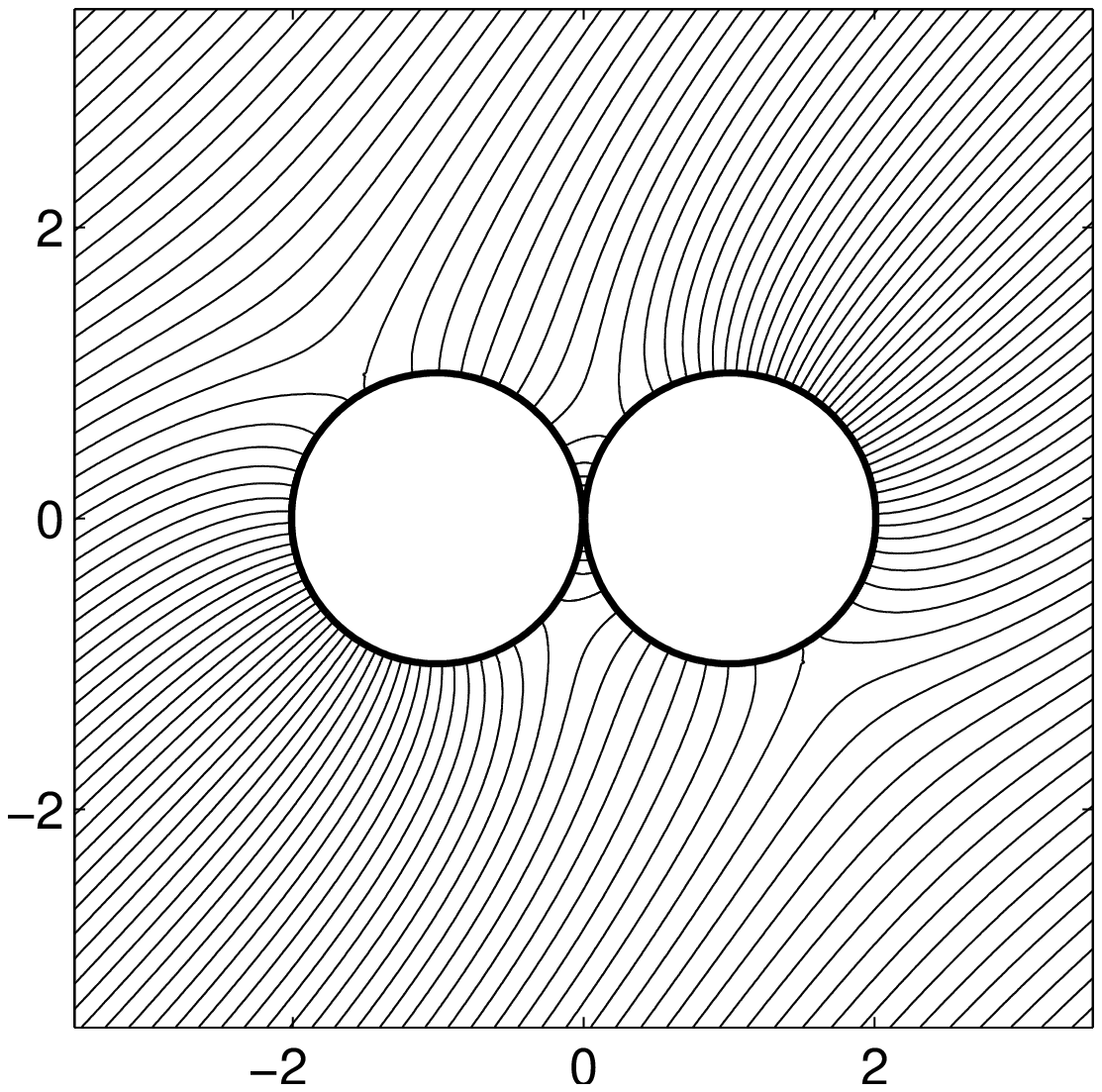, height=4cm}\hskip 1cm
\epsfig{figure=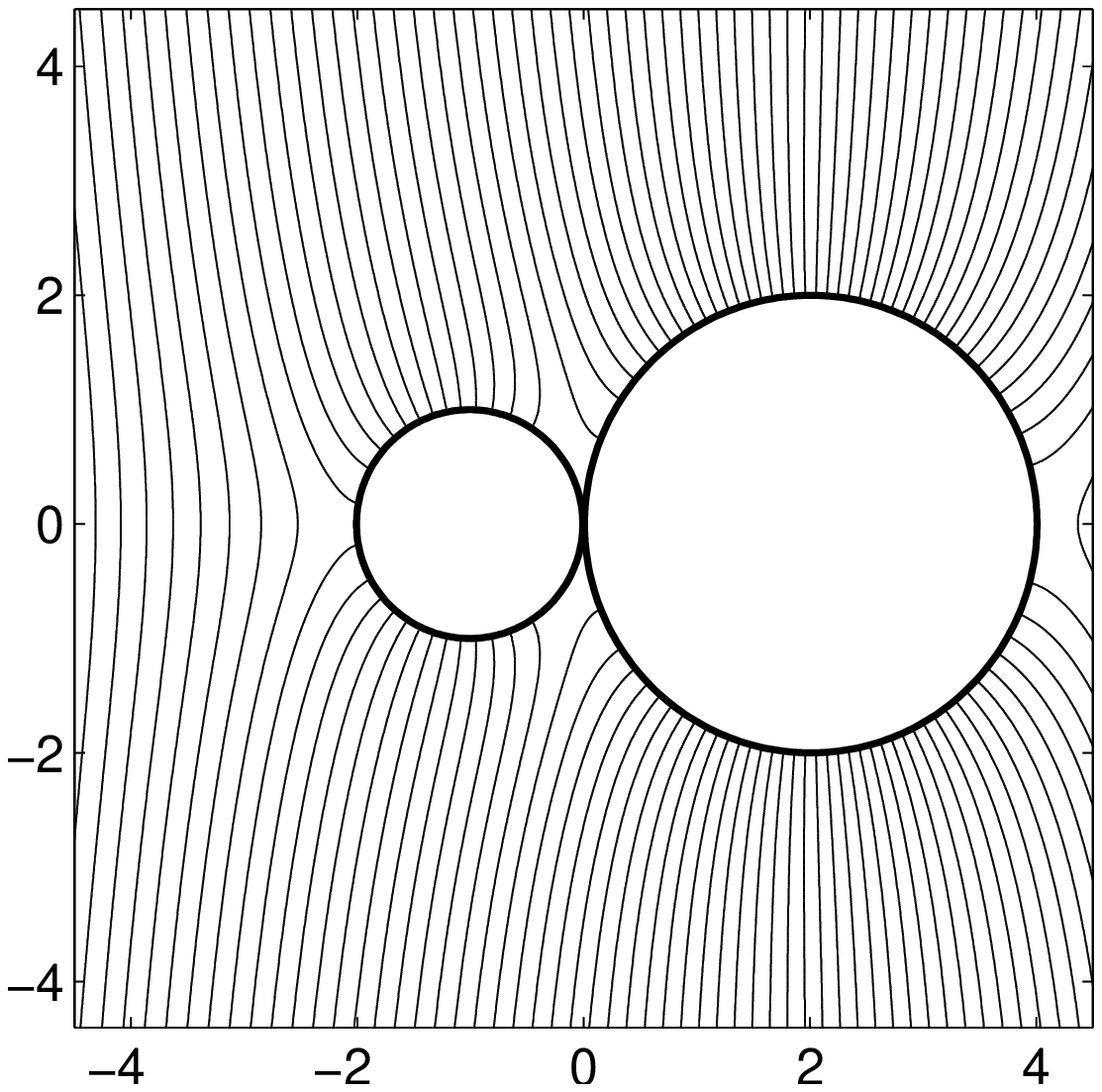, height=4cm}\\[.5cm]
\end{center}
\caption{Level curves of the free space insulating case. The entire harmonic function $H(\Bx) = x_1$ in the upper-left and the lower-right figure, and $H(\Bx) = x_2$ in the upper-right one, and
$H(\Bx) =x_1 -x_2$ in the lower-left one.}\label{level_curve_insul}
\end{figure}

In Figure \ref{level_curve} and \ref{level_curve_insul}, the uniformly spaced contour level curves are drawn for the free space conducting and insulating case, respectively.
The distance $\epsilon=0.0156$ and the number of grid points on each disk is 256.
The radii are $r_1=r_2=1$ except the lower-right figure where $r_2=2$.
The entire harmonic function $H(\Bx) = x_1$ in the upper-left and the lower-right figure, and $H(\Bx) = x_2$ in the upper-right one, and
$H(\Bx) =x_1 -x_2$ in the lower-left one.

\begin{thebibliography}{00}

%\bibitem{A}{ L. Ahlfors}, {Complex Analysis}, Third ed., McGraw-Hill, New
%York, 1979.

\bibitem{ADKL} { H. Ammari, G. Dassios, H. Kang H and  M. Lim}, { Estimates for the electric field in the presence of adjacent perfectly conducting spheres}, {Quat. Appl. Math.} { 65} (2007), 339--355.

\bibitem{book2} {H. Ammari and H. Kang}, {Polarization and moment
tensors with applications to inverse problems and effective medium
theory}, Applied Mathematical Sciences, Vol. 162, Springer-Verlag,
New York, 2007.

\bibitem{AKL} { H. Ammari, H. Kang and M. Lim}, {Gradient estimates for solutions to the conductivity problem}, Math. Ann. 332(2) (2005),~277--286.

\bibitem{AKLLL} { H. Ammari, H. Kang, H. Lee, J. Lee and M. Lim}, { Optimal bounds on the gradient of solutions to conductivity problems}, J. Math. Pures Appl. 88 (2007),~307--324.

\bibitem{AKLLZ} H. Ammari, H. Kang, H. Lee, M. Lim and H. Zribi, Decomposition theorems and fine estimates for electrical fields in the presence of closely located circular inclusions, Jour. Diff. Equa., 247 (2009), 2897-2912.

%\bibitem{Bao_thesis}{ E. S. Bao}, { Gradient estimates for the conductivity problems}, thesis, Rutgers University.

\bibitem{bab} I. Babu\u{s}ka, B. Andersson, P. Smith and K. Levin, Damage analysis of fiber composites. I. Statistical
analysis on fiber scale, Comput. Methods Appl. Mech. Engrg. 172 (1999), 27--77.

\bibitem{BC}{  B. Budiansky and G. F. Carrier}, { High shear stresses in stiff fiber composites}, Jour. Appl. Mech. 51 (1984),  733-735.

\bibitem{BLY} { E. S. Bao, Y.Y. Li and B. Yin}, { Gradient estimates for the perfect conductivity problem}, Arch. Rat. Mech. Anal. 193 (2009), 195-226.

\bibitem{BLY2} {E. S. Bao, Y.Y. Li and B. Yin}, { Gradient Estimates for the perfect and insulated conductivity problems with multiple inclusions}, Comm.  Part. Diff. Equa. 35 (2010), 1982--2006.

\bibitem{BV} { E. Bonnetier and  M. Vogelius}, {An elliptic regularity result for a composite medium with ``touching" fibers of circular cross-section}, SIAM Jour. Math. Anal. 31 No 3 (2000), 651--677.

\bibitem{CG} H. Cheng and L. Greengard,
\newblock A method of images for the evaluation of electrostatic
fields in systems of closely spaced conducting cylinders,
\newblock SIAM J. Appl. Math. 58 (1998), 122-141.

\bibitem{Greengard} L. Greengard and M. Moura,
\newblock On the numerical evaluation of electrostatic fields in
composite materials, \newblock Acta Numerica (1994), 379--410.

\bibitem{FKS99} { E. Fabes, H. Kang, and J.K. Seo}, { Inverse conductivity problem
with one measurement: Error estimates and approximate
identification for perturbed disks},
 SIAM J. Math. Anal., 30 (1999), 699--720.

\bibitem{KS96} {H. Kang and J.K. Seo}, {Layer potential technique for the inverse
conductivity problem}, Inverse Problems, 12 (1996), 267--278.

\bibitem{KS2000}
{H. Kang and J.K. Seo}, {Recent progress in the inverse conductivity problem with
single measurement}, in Inverse Problems and Related Fields, CRC
Press, Boca Raton, FL, (2000), 69--80.
\bibitem{K} { J.B. Keller}, {Stresses in narrow regions}, Trans. ASME J. Appl. Mech. 60 (1993), 1054--1056.

\bibitem{mar} X. Markenscoff, Stress amplification in vanishing
small geometries, Computational Mechanics 19 (1996), 77--83.

\bibitem{LN} {Y.Y. Li and L. Nirenberg}, { Estimates for elliptic system from composite material}, Comm. Pure Appl. Math., LVI (2003), 892--925.


\bibitem{LV} { Y.Y. Li and M. Vogelius}, {Gradient estimates for solution to divergence form elliptic equation with discontinuous coefficients}, Arch. Rat. Mech. Anal. 153 (2000), 91--151.

\bibitem{LY} {M. Lim and K. Yun}, {Blow-up of electric fields between closely spaced spherical perfect conductors}, Comm.  Part. Diff. Equa. 34 (2009), 1287--1315.

\bibitem{LY2} {M. Lim and K. Yun}, {Strong influence of a small fiber on shear stress in fiber-reinforced composites}, Jour. Diff. Equa. 250 (2011), 2402--2439.

\bibitem{Y} {K. Yun},  {Estimates for electric fields blown up between closely adjacent conductors with arbitrary shape}, { SIAM Jour. Appl. Math.} 67 No 3 (2007), 714--730.

\bibitem{Y2} {K. Yun}, {Optimal bound on high stresses occurring between stiff fibers with arbitrary shaped cross sections}, Jour. Math. Anal. Appl. 350 (2009), 306-312.



\end{thebibliography}
\end{document}